\documentclass[journal]{IEEEtran}      %
\usepackage{amsmath,amsfonts,amssymb}
\usepackage{algorithmic}
\usepackage{algorithm}
\usepackage{array}
\usepackage[caption=false,font=normalsize,labelfont=sf,textfont=sf]{subfig}
\usepackage{textcomp}
\usepackage{stfloats}
\usepackage{url}
\usepackage{verbatim}
\usepackage{graphicx}
\usepackage{cite}

\usepackage{epsfig} %
\usepackage{physics}
 \usepackage{multirow}
 \usepackage{epstopdf}
 \usepackage{tikz}
\usetikzlibrary{arrows.meta, %
                bending,     %
                patterns     %
               }

\title{\LARGE \bf
Speeding up VSLMS adaptation algorithms using dynamic adaptation gain: Analysis and Applications *
}

\author{Ioan Doré Landau$^{a}$, Dariusz Bismor$^{b}$, Tudor-Bogdan Airimitoaie$^{c}$, Bernard Vau$^{d}$, Gabriel Buche$^{a}$%
\thanks{*This work was not supported by any organization}%
\thanks{$^{a}$Ioan Doré Landau and Gabriel Buche are with the Univ. Grenoble Alpes, CNRS, Grenoble INP, GIPSA-lab, 38000 Grenoble, France
        {\tt\small firstname.lastname@gipsa-lab.grenoble-inp.fr}}%
\thanks{$^{b}$Dariusz Bismor is with the Silesian University of Technology, ul. Akademicka 2A, 44-100 Gliwice, Poland
                {\tt\small dariusz.bismor@polsl.pl}}%
\thanks{$^{c}$Tudor-Bogdan Airimitoaie is with the Univ. Bordeaux, CNRS, Bordeaux INP, IMS, 33405 Talence, France
        {\tt\small tudor-bogdan.airimitoaie@u-bordeaux.fr}}%
\thanks{$^{d}$Bernard Vau is with Exail, 12 avenue des Coquelicots, 94385 Bonneuil sur Marne, France
        {\tt\small bernard.vau@exail.com}}%
}

\newtheorem{thm}{Theorem}[section]

\newtheorem{lem}[thm]{Lemma}

\begin{document}

\maketitle
\thispagestyle{empty}
\pagestyle{empty}

\begin{abstract}
The paper explores the use of dynamic adaptation gain/step size (DAG) for improving the adaptation  transient performance of variable step-size LMS (VS-LMS) adaptation algorithms. A generic form for the implementation of the DAG  within the VS-LMS algorithms is provided. The properties of the VS-LMS algorithms using dynamic adaptation gain are discussed in detail. Stability issues in deterministic environment and convergence properties in stochastic environment are examined. A transient performance analysis is proposed.  Criteria for  the selection of the coefficients of the DAG filter are provided.
The potential of the VS-LMS adaptation algorithms using a DAG is then illustrated by simulation results (adaptive line enhancer, filter identification)  and experimental results obtained on a relevant adaptive active noise attenuation system. 
\end{abstract}

\section{INTRODUCTION}
The modern development of adaptation techniques in automatic control and signal processing started at the end of the fifties and beginning of the sixties (20th century). The paper \cite{Widrow60} introduced a gradient-based adaption algorithm in the discrete-time, later named the least mean squares (LMS). While the choice of the adaptation gain/step size for assuring the stability of the system was an open problem, interesting applications in the field of signal processing have been done. The paper \cite{Widrow75} gives an account of the applications of the LMS algorithm up to 1975.

In automatic control, the first attempt to synthesize adaptation algorithms has been probably the paper \cite{Whitaker58}, where a continuous-time formulation of a gradient type algorithm has been proposed. Unfortunately, dealing with feedback control systems to which a non-linear/time-varying loop (the adaptation loop) is added, raised crucial stability issues. The problem of the choice of the adaptation gain (step size) assuring the stability of the full system is fundamental. Therefore, the research paradigm in control was directed toward synthesis of adaptive algorithms guaranteeing the stability of the full system for any (positive) value of the adaptation gain/step size. 
 Discrete time adaptation  algorithms assuring global asymptotic stability for any values of the adaptation gain were available since 1971 (\cite{LandauAsilomar71,LandauIFAC73_1,LandauIFAC75}). The concepts of "a priori" and "a posteriori" adaptation error emerged as key points for understanding the stability issues in the discrete-time context. Use of this type of algorithms in signal processing have been reported in \cite{Johnson79}, \cite{Larimore80} among other references.
   These algorithms, derived from stability considerations, can be interpreted in the scalar case as gradient type algorithms trying to minimize a quadratic criterion in terms of the "a posteriori" prediction error (\cite{LandauSpringer11,Landau90c}).

The signal processing community has concentrated its efforts in developing variable step size/adaptation gain algorithms in a scalar context (more exactly using a diagonal matrix adaptation gain) in order to improve the performance of the LMS algorithm. An exhaustive account of the various variable step-size LMS (VS-LMS) algorithms is provided by \cite{Bismor16}, where a unified presentation is done as well as an extensive comparison of the algorithms by applying them to a number of significant problems. It turns out that the adaptation algorithms developed in automatic control from a stability point of view can also be interpreted as "variable step size" LMS algorithms.

When using adaptive/learning recursive algorithms there is an important problem to be addressed: the compromise between alertness (with respect to environment changes -- like plant or disturbance characteristics) and stationary performances when using a constant value for the adaptation gain/step size. Accelerating the adaptation transient without augmenting the value of the adaptation gain/step size is a challenging problem.

Recently, the concept of dynamic adaptation gain (DAG) has been introduced in  \cite{LandauCDC22,LandauAuto23,LandauJSV23} as a way to accelerate significantly the adaptation transients without modifying the steady state (asymptotic) properties of an algorithm for a given adaptation gain/step size. The correcting term in the adaptation algorithm is first filtered before its use for the estimation of a new parameter value. With an appropriate choice of the parameters of this filter, which should be characterized by a strictly positive real (SPR) transfer function, a significant improvement of the adaptation transient is obtained. The design of this filter is well understood and design tools are available.

The main objective of this paper is to show that the DAG introduced in the context of stability based adaptation algorithms can be successfully applied to VS-LMS adaptation algorithms leading to similar significant acceleration of the adaptation transients. This will be based on a theoretical analysis and will be illustrated by simulations (adaptive line enhancer, filter identification) and real-time experiments on an adaptive feedforward noise attenuator (a silencer).

The paper is organized as follows: Section \ref{DAG} will introduce the concept of dynamic adaptation gain in the context of VS-LMS adaptation algorithms. Section \ref{design} will focus on the analysis and design of the dynamic adaptation gain filter.
 Section \ref{Sim} will present simulations results for an adaptive line enhancer and for filter identification. Section \ref{Exp} will present the experimental results obtained with three algorithms on an adaptive feedforward noise attenuator.
\section{Introducing the DAG-VS-LMS algorithms}\label{DAG}
\begin{figure}[!htb]
    \begin{center}
    \includegraphics[width=0.6\columnwidth]{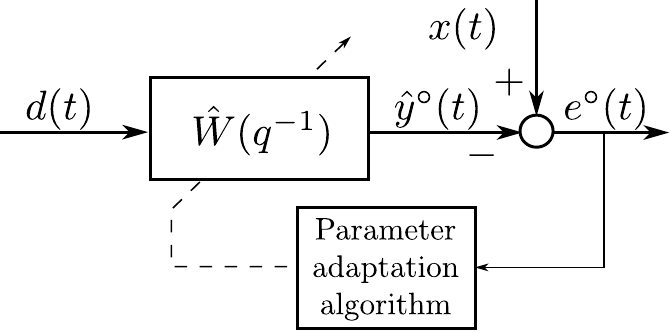}
    \caption{Least mean squares (LMS) adaptive filtering problem.}
    \label{fig_LMS_basic_scheme}
    \end{center}
\end{figure}

Variable step size least mean squares (VS-LMS) algorithms are very popular in the field of signal processing and the field of active vibration and noise control. There is a strong similarity with some of the algorithms used in the adaptive control and recursive system identification. The VS-LMS algorithms (which are improvements of the original LMS algorithm) will be briefly reviewed in order to add the dynamic adaptation gain/step size introduced in \cite{LandauCDC22}. The aim of the LMS parameter adaptation/learning algorithm (PALA) is to drive the parameters of an adjustable model in order to  minimize a quadratic criterion in terms of the prediction error (difference between real data and the output of the model used for prediction).

The basic block diagram illustrating the LMS algorithm's operation is shown in Fig.~\ref{fig_LMS_basic_scheme}. The adaptive filter $\hat{W}(q^{-1})$ is fed with the input sequence\footnote{Variables in discrete-time are denoted as $s(t)$, where $t$ gives the integer number of sampling periods and  it is related to the continuous-time by $\tau=t\cdot T_s$, where $T_s$ is the sampling period.} $d(t)$. The output of the filter, $\hat{y}^\circ(t)$, is compared with the desired signal, $x(t)$, to compute the error signal $e^\circ(t)$. The LMS algorithm adjusts the weights of the $\hat{W}(q^{-1})$ filter to minimize the error. 

Consider that the desired signal can be described by:
\begin{align}
  x(t)=\mathbf{w}^T\mathbf{r}(t),
  \label{eq:output}
\end{align}
where the \textit{parameter vector} and the \textit{measurement vector} are denoted by
\begin{align}
  \mathbf{w}^T&=\left[ w_0,w_1,\ldots,w_{n_W} \right],\text{ and}\\
  \mathbf{r}^T(t)&=\left[ d(t),d(t-1),\ldots,d(t-n_W) \right],
\end{align}
respectively. The adjustable prediction model of the adaptive filter will be described by:
\begin{align}
  \hat{y}^\circ(t)=\hat{\mathbf{w}}^T(t-1)\mathbf{r}(t),
\end{align}
where $\hat{y}^\circ(t)$ is termed the \textit{a priori} predicted output depending upon the values of the estimated parameter vector $\mathbf{w}$ at instant $t-1$:
\begin{align}
  \hat{\mathbf{w}}^T(t-1)&=\left[ \hat{w}_0(t-1),\hat{w}_1(t-1),\ldots,\hat{w}_{n_W}(t-1) \right].
\end{align}
It is very useful to consider also the \emph{a posteriori} predicted output computed on the basis of the new estimated parameter vector at $t$, $\hat{\mathbf{w}} (t)$, which will be available somewhere between $t$ and $t+1$. The \emph{a posteriori} predicted output will be given by:
\begin{align}
\hat{y}(t)&=\hat{\mathbf{w}}^T(t) \mathbf{r}(t)
\label{a postout}
\end{align}
One defines an \emph{a priori} prediction error as:
\begin{equation}
 e^\circ(t)=x(t)-\hat{y}^\circ(t)= [\mathbf{w} - \hat{\mathbf{w}}(t-1)]^T\mathbf{r}(t) \label{eq_apriori}
\end{equation}
and an \emph{a posteriori} prediction error\index{a posteriori prediction error} as:
\begin{equation}
e(t)=x(t)-\hat{y}(t)= [\mathbf{w} - \hat{\mathbf{w}}(t)]^T\mathbf{r}(t) \label{eq_aposteriori}
\end{equation}
The VS-LMS algorithms update the filter taps according to the formula (see \cite{Bismor16}):
\begin{align}
  \hat{\mathbf{w}}(t)=\hat{\mathbf{w}}(t-1)+\mu(t) \mathbf{r}(t) e^\circ(t),
  \label{eq_vslms}
\end{align}
where $\mu(t)$ is the variable step-size parameter. In the standard form of the LMS algorithm, the step-size has a constant value $\mu(t)=\mu$. Large values of $\mu$ allow for fast adaptation, but also give large excess mean square error (EMSE, see \cite{Bismor16}). Too large step-sizes may lead to the loss of stability (see \cite {LandauSpringer11} for a discussion). On the other hand, too small step-sizes give slow convergence, which in many practical applications is not desirable. A variable step-size $\mu(t)$ can provide a compromise. The first VS-LMS algorithm was the Normalized LMS (NLMS) algorithm proposed in 1967 independently by \cite{Albert67,Nagumo67}, which uses the following equation for the step-size:
\begin{align}
  \mu(t)=\frac{\mu}{\delta+\mathbf{r}(t)^T\mathbf{r}(t)},
\end{align}
where $\delta$ is a very small scalar value used in order to avoid division by zero (typically of the order of $10^{-16}$).
\footnote{The case of NLMS using larger values of $\delta$ (for example $\delta=1$) is discussed in \cite{Goodwin84a}.}

In adaptive control, where the stability of the full adaptive control system is considered as a fundamental issue, VS-LMS algorithms have been developed from a stability point of view. However, when scalar type adaptation gain is used (i.e. diagonal matrix adaptation gain), these algorithms  can be interpreted as gradient type algorithms for the minimization of a quadratic error in terms of the \textit{a posteriori} prediction error defined in \eqref{eq_aposteriori} - see \cite{LandauSpringer11} for details. See also \cite{LandauJSV23} for a derivation of this type of algorithm and its application to active noise control. This algorithm will be termed PLMS \footnote{In \cite{Johnson79} a close algorithm is termed HARF.} (to distinguish it with respect to the LMS and the other VS-LMS that use the \textit{a priori }prediction error). The algorithm has the form:
\begin{align}
  \hat{\mathbf{w}}(t)&=\hat{\mathbf{w}}(t-1)+\mu \mathbf{r}(t) e(t) \label{eq:plms1}\\
  &=\hat{\mathbf{w}}(t-1)+\mu(t) \mathbf{r}(t) e^\circ(t),
  \label{eq_plms}
\end{align}
where:
 \begin{align}
  \mu(t)=\frac{\mu}{1+\mu\mathbf{r}(t)^T\mathbf{r}(t)},
  \label{eq_mu}
\end{align}
$e(t)$ can be computed from $e^\circ(t)$. One has using Eq. \eqref{eq:plms1}:
\begin{equation}
e(t)= [\mathbf{w} - \hat{\mathbf{w}}(t)]^T\mathbf{r}(t)= e^\circ(t)-\mu\mathbf{r}^T\mathbf{r}(t) e(t)
\end{equation}
yielding: 
\begin{equation}
e(t)=\frac{e^\circ(t)}{1+\mu\mathbf{r}(t)^T\mathbf{r}(t)}
\label{eq:e-e0}
\end{equation}

When using the \textit{dynamic adaptation gain/step size} (DAG),  Eq. \eqref{eq_vslms} of the VS-LMS algorithms will take the form:
\begin{equation}
\hat{\mathbf{w}}(t)=\hat{\mathbf{w}}(t-1)+\frac{C(q^{-1})}{D'(q^{-1})}\left[ \mu(t) \mathbf{r}(t) e^\circ(t) \right]
\label{eq:dynamic1}
\end{equation}
where\footnote{The complex variable $z^{-1}$ will be used for characterizing the system's behaviour in the frequency domain and the delay operator $q^{-1}$ will be used for describing the system's behavior in the time domain.} $\frac{C(q^{-1})}{D'(q^{-1})}$ is termed the ``dynamic adaptation gain/step size" (DAG) and has the form:
\begin{equation}\label{eq_theta_DAG}
H_{DAG}(q^{-1})=\frac{C(q^{-1})}{D'(q^{-1})}=\frac{1+c_{1}q^{-1} + \ldots+c_{n_C}q^{-n_C}}{1-d'_{1}q^{-1}-\ldots -d'_{n_{D'}}q^{-n_{D'}}}
\end{equation}
The effective implementation of the algorithm given in Eq. \eqref{eq:dynamic1} leads to:
\begin{multline}
 \hat{\mathbf{w}} (t) =d_1\hat{\mathbf{w}}(n-1)+ d_2\hat{\mathbf{w}}(n-2)+\ldots+d_{n_D}\hat{\mathbf{w}}(n-n_D) \\
 +\mu(t)\mathbf{r}(t)e^\circ(t)+ c_1\mu(n-1)\mathbf{r}(n-1)e^\circ(n-1)+ \\
 + \ldots +c_{n_C}\mu(n-n_C)\mathbf{r}(n-n_C)e^\circ(n-n_C)
 \label{eq:ARMA3}
\end{multline}
where ($n_D=n_{D'}+1$):
\begin{equation}
d_i=(d'_i-d'_{i-1})~~; i=1,... n_D; d'_0=-1,~d'_{n_D}=0
\label{ARIMA2}
\end{equation}
To implement the algorithm, one needs a computable expression for $e^\circ(t)$. This is obtained by computing $\hat{y}^\circ(t)$ in \eqref{eq_apriori} as\footnote{For low orders $n_C$ and $n_{D'}$, $\hat{\mathbf{w}_0}(t-1)$ can efficiently be approximated by $\hat{\mathbf{w}}(t-1)$.} 
\begin{align}
  \hat{y}^\circ(t)=\hat{\mathbf{w}}_0^T(t-1)\mathbf{r}(t),
  \label{eq_z0DAG}
\end{align}
where
\begin{multline}
  \hat{\mathbf{w}}_0^T(t-1)=d_1\hat{\mathbf{w}}(t-1)+ d_2\hat{\mathbf{w}}(t-2)+\ldots+d_{n_D}\hat{\mathbf{w}}(t-n_D) \\
  +c_1\mu(t-1)\mathbf{r}(t-1)e^\circ(t-1)+ \\
  + \ldots +c_{n_C}\mu(t-n_C)\mathbf{r}(t-n_C)e^\circ(t-n_C).
\end{multline}

\subsection*{Relations with other algorithms}
Many algorithms have been proposed for accelerating the speed of convergence of the adaptation algorithms derived using the "gradient rule". The algorithm of \eqref{eq:dynamic1} is termed ARIMA (Autoregressive with Integrator Moving Average). As discussed in \cite[Section 8]{LandauAuto23}, a number of well known algorithms are particular cases of the ARIMA algorithm. The various algorithms described in the literature are of MAI (Moving Average with Integrator) form or ARI (Autoregressive with Integrator) form. The MAI form includes ``Integral+ Proportional" algorithm \cite{LandauSpringer11,AirimitoaieAuto13} ($c_1\neq  0, c_i=0,\forall~i>1$, $d'_i=0,\forall~i>0$),
``Averaged gradient" ($c_i, i=1,2,...$, $d'_i=0,\forall~i>0$) \cite{Schmidt18,Pouyanfar17}. The ARI form includes ``Conjugate gradient" and ``Nesterov" algorithms \cite{Livieris14,AnnaswamyCDC19} ($c_i=0, i=1,2,..,d'_1 \neq  0, d'_i=0,i>1)$ as well as the ``Momentum back propagation" algorithm \cite{Jacobs88} which corresponds to the conjugate gradient plus a normalization of $\mu$ by $(1-d'_1)$\footnote{There are very few indications how to choose the various coefficients in the above mentioned algorithms.}.
A particular form of the ARIMA algorithms termed ``ARIMA2" ($c_1,c_2 \neq 0, c_i=0,\forall~i>2,d'_1\neq 0,d'_i=0,\forall~i>1$) will be studied subsequently and evaluated experimentally.\footnote{The algorithms mentioned above can be viewed as particular cases of the ARIMA2 algorithm.}
\section{Analysis and design of the dynamic adaptation gain/step size}\label{design}
\subsection{Performance issues}
The dynamic adaptation gain/learning rate will introduce a frequency-dependent phase distortion on the gradient.  Assume that the algorithms should operate for all frequencies in the range: $0$ to $0.5f_s$ ($f_s$ is the sampling frequency). 
Assume that the gradient of the criterion to be minimized contains a single frequency.
 In order to minimize the criterion, the phase distortion introduced by the dynamic adaptation gain/learning rate should be less than $90^\circ $ at all the frequencies. 
In other terms, the transfer function $\frac{C(z^{-1})}{D'(z^{-1})}$ should be strictly positive real (SPR).
 In order that a transfer function be strictly positive real, it must first have its zeros and poles inside the unit circle. One has the following property:
\begin{lem}\label{lem_I_equal_0}
  Assume that the polynomials $C(z^{-1})$ and $D'(z^{-1})$ have all their zeros inside the unit circle, then:
  \begin{equation}
  I=\int_{0}^{\pi}\log\left( \left|\frac{C(e^{-i\omega})}{D^{'}(e^{-i\omega} )}\right | \right)\dd  \omega =0.
  \label{eq:integral}
  \end{equation}
\end{lem}
The proof relies on the Cauchy Integral formula (see \cite{LandauAuto23}).

This result allows to conclude that the average gain of a SPR filter over the frequency range $0$ to $0.5f_s$ is 0~dB, i.e. on the average this filter will not modify the adaptation gain/step size. It is just a frequency weighting of the adaptation gain/step size. A it can be seen in Fig. \ref{fig:dagbode},  Fig. \ref{fig:C_D'} (Section \ref{Sim}) and in Fig.~\ref{fig_bode_comp} (Section \ref{Exp}), these filters are SPR low band pass filters with an average gain of $0$ dB over the frequency range $0$ to $0.5 f_s$.
 This means that if the frequency content of the gradient is in the low frequency range, the effective adaptation gain/learning rate will be larger than $\mu$, which should have a positive effect upon the adaptation/learning transient. In particular, the DAG which has the largest gain in low frequencies  should provide the best performance  (This is indeed the case---see Sections \ref{Sim} and \ref{Exp}).

Since we need to have a DAG which is SPR,  we will provide subsequently the tools for design of a SPR DAG. We will consider the case of the ARIMA2 algorithm introduced in \cite{LandauCDC22}. The DAG in this case will have the form:
\begin{align}\label{eq_H_DAG}
H_{DAG}=\frac{C(q^{-1})}{D'(q^{-1})}=\frac{1+c_1q^{-1}+c_2q^{-2}}{1-d'_1 q^{-1}}
\end{align}
A criterion for the selection of $c_1$, $c_2$ and $d'_1$ in order that the DAG be SPR \cite{LandauAuto23} is given in Appendix \ref{Appen_SPR}.

From the conditions of Lemma \ref{lemma_SPR} (Appendix \ref{Appen_SPR}), closed contours in the plane $c_2-c_1$ can be defined for the different values of $d'_1$ allowing to select $c_1$ and $c_2$ for a given value of $d'_1$ such that the DAG be SPR. 
Note that a necessary condition for the selection of the parameters $c_1, c_2, d'_1$ is that both the denominator and the numerator of the $H_{DAG}$ should be asymptotically stable.
\subsection{Stability issues for unconstrained values of the adaptation gain $\mu>0$}
If one wants to use VS-LMS algorithms with large values of the adaptation gain $\mu$, the stability analysis of the resulting scheme using a dynamic adaptation gain is an important issue. For the case of the PLMS algorthm, this analysis has been carried in detail in \cite{LandauAuto23, LandauJSV23}\footnote{For $\mu=1$ and $\delta=1$, PLMS and NLMS are identical}. For the convenience of the reader, we will indicate the main result subsequently.\\
Equation \eqref{eq:dynamic1} can be expressed also as:
\begin{equation}
(1-q^{-1})\hat{\mathbf{w}}(t+1)=+\frac{C(q^{-1})}{D'(q^{-1})}\left[ \mu(t) \mathbf{r}(t) e^\circ(t+1) \right] \nonumber
\label{eq:dynamic1''}
\end{equation}
leading to:
\begin{equation}
\hat{\mathbf{w}}(t+1)=H_{PAA}(q^{-1})[\mu\mathbf{r}(t) e(t+1)],
\label{eq:genalg3}
\end{equation}
where  $e(t+1)=e^\circ(t+1)(1+\mu\phi^T(t)\phi(t))^{-1}$ and $H_{PAA}=(1-q^{-1})^{-1}H_{DAG}$ is a MIMO diagonal transfer operator having identical terms. All the diagonal terms are described by:
\begin{align}
H_{PAA}^{ii}(q^{-1})&=\frac{1+c_{1}q^{-1} +\ldots+c_{n_C}q^{-n_C}}{(1-q^{-1})(1-d'_{1}q^{-1}-\ldots-d'_{n_{D'}}q^{-n_D'})}\nonumber \\
&=\frac{C(q^{-1})}{(1-q^{-1})D'(q^{-1})}=\frac{C(q^{-1})}{D(q^{-1})}.
\label{eq:ARIMA1}
\end{align}
The relation between the coefficients of polynomials $D$ and $D'$ is given  in \eqref{ARIMA2}.

One has the following result:
\begin{thm}\label{thm_PR}
  For the system described by Equations~\eqref{eq:output} through \eqref{eq_aposteriori} using the PLMS  algorithm of \eqref{eq:dynamic1} and \eqref{eq_theta_DAG}
  one has $\lim_{t \to \infty} e(t+1) =0$ for any positive adaptation gain $\mu$ and any initial conditions $\hat{\mathbf{w}} (0), e (0)$, if $H_{PAA}^{ii}(z^{-1})$  given in \eqref{eq:ARIMA1} is a positive real (PR) transfer function with a pole at $z=1$.
\end{thm}

The proof is given in \cite{LandauAuto23}. %

For the particular case of the ARIMA2 algorithm, the coefficients $c_1, c_2$ and $d'_1$ should be chosen such that the DAG is SPR and the $H^{ii}_{PAA}$ is positive real (PR), i.e.
\begin{align}
H_{PAA}^{ii}&=\frac{1+c_1q^{-1}+c_2q^{-2}}{1-d_1q^{-1}-d_2 q^{-2}}=\frac{1+c_1q^{-1}+c_2q^{-2}}{(1-q^{-1})(1-d'_1 q^{-1})}
\label{eq:ARIMA7}
\end{align}
should be PR. A criterion for the selection of the coefficients $c_1, c_2$ and $d_{1}'$ can be found in \cite{LandauAuto23}. Closed contours cans be defined in the $c_2-c_1$ plane for given values of $d_{1}'$. 

The adaptive/learning system considered in the Theorem~\ref{thm_PR}, leads to an equivalent feedback representation where the equivalent feedforward path is a constant positive gain and the equivalent feedback path features the $H_{PAA}$ (see \cite{LandauAuto23}). The feedback path is passive under the condition that $H^{ii}_{PAA}$ is positive real and the feedforward pass is strictly passive allowing to conclude upon the stability of the closed loop. However, in a number of  cases (like output error configurations) the equivalent feedforward path may be a transfer operator.  In such situations in addition to the  PR condition upon the $H_{PAA}$, there will be an additional SPR condition upon the transfer operator characterizing the equivalent feedforward path.

\begin{figure}[htb!]%
    \begin{center}
    \includegraphics[width=\columnwidth]{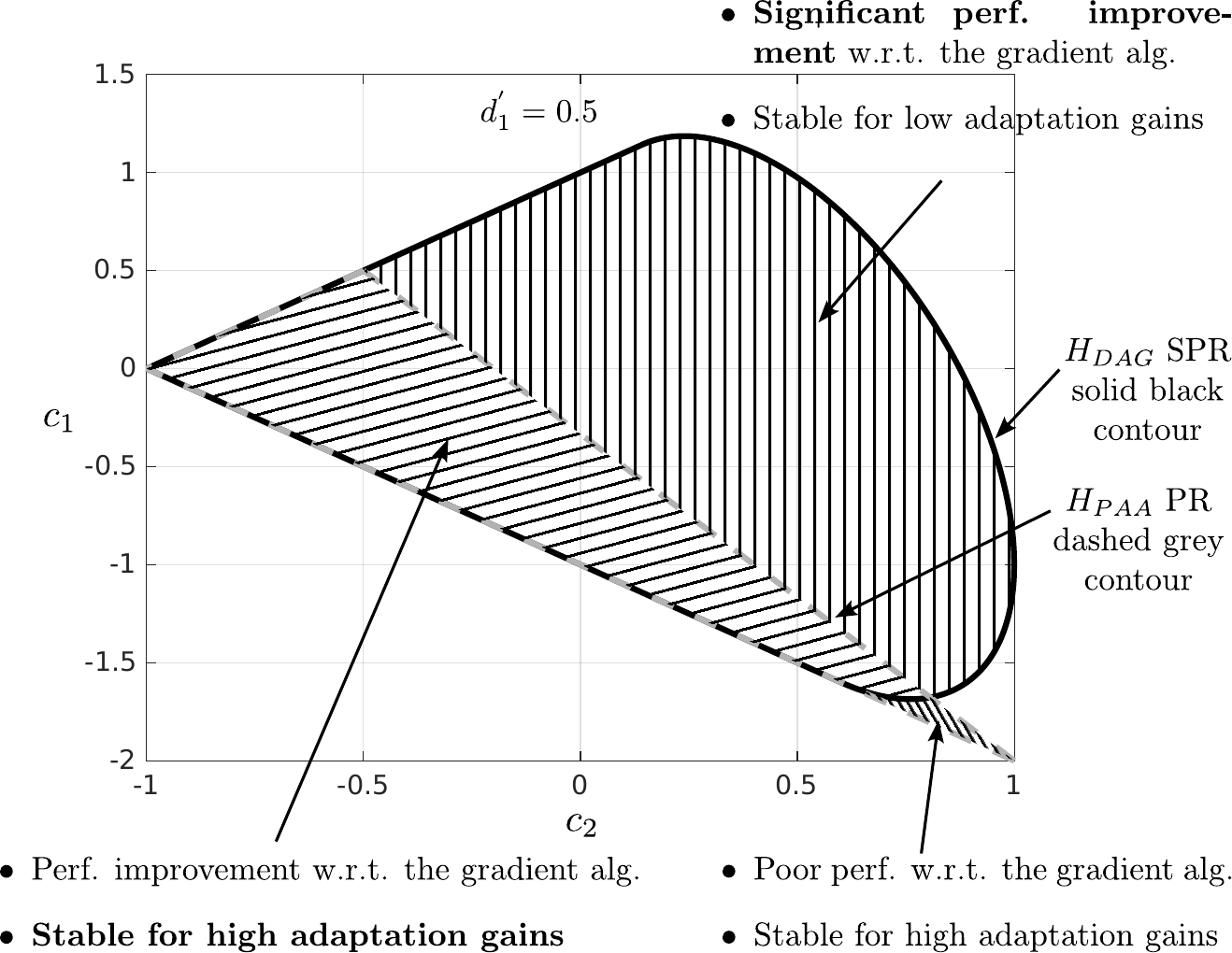}
    \caption{Intersection in the plane $c_1 - c_2$ of the contour  $H_{PAA}=PR$ with the contour  $H_{DAG}=SPR$ for $d'_1=0.5$.}
    \label{fig_GADww}
    \end{center}
\end{figure}

It is interesting to see intersections of the contours ensuring the SPR of the $H^{ii}_{DAG}$ with the contours assuring that $H^{ii}_{PAA}$ is PR. Such an intersection is shown in Fig.\ref{fig_GADww}. From this figure one can conclude that not all the SPR $H_{DAG}$ will lead to a PR $H_{PAA}$. In such cases, the performance is improved for low adaptation gains, but asymptotic stability cannot be guaranteed for high adaptation gain values.    Fig.~\ref{fig_GADww} shows also that there is a region where despite  that $H_{PAA}$ is PR, $H_{DAG}$ is not SPR. For such configurations, large adaptation gains can be used but the adaptation transient is slower than for the basic gradient algorithm.

\subsection{Stability issues for low values of the adaptation gain/step size $\mu>0$}

For small values of the adaptation gains/learning rates, the passivity/stability condition can be relaxed using \textit{averaging}  \cite{Anderson86}.
 Using the results of \cite{LandauAuto11}, under the hypothesis of an input signal spanning all the frequencies up to half of the sampling frequency, \textit{passivity in the average} will be assured if the frequency interval where $H^{ii}_{PAA}$ is not positive real is smaller than the frequency interval where $H^{ii}_{PAA}$ is positive real. In fact, the most important is that the $H^{ii}_{PAA}$ is PR in the frequency region of operation (mainly defined by the spectrum of the input signals to the systems). However, more specific results can be obtained using \textit{averaging} as shown next\footnote{The subsequent analysis is valid for most of the VS-LMS algorithms.}.\\
 Defining the parameter error as:
 \begin{equation}
 \tilde{\mathbf{w}}(t)= \hat{\mathbf{w}}(t)- \mathbf{w}
 \end{equation}
 from Eq. \eqref{eq:genalg3} on gets:
 \begin{equation}
\tilde{\mathbf{w}}(t+1)=\tilde{\mathbf{w}}(t)+FH_{DAG}(q^{-1})[\mathbf{r}(t) e(t+1)];~~F=\mu I
\label{eq:par_err3}
\end{equation}
Taking into account that :
\begin{equation}
e(t+1)=\mathbf{r}^T[\mathbf{w}-\hat{\mathbf{w}}(t+1)]=-\mathbf{r}^T\tilde{\mathbf{w}}(t+1)
\end{equation}
 Eq. \eqref{eq:par_err3} becomes:
 \begin{align}
\tilde{\mathbf{w}}(t+1)
&=\tilde{\mathbf{w}}(t)-FH_{DAG}(q^{-1})[\mathbf{r}(t)\mathbf{r}^T\tilde{\mathbf{w}}(t+1)]
\label{eq:par_err4}
\end{align}

\begin{figure}[htb!]%
    \begin{center}
    \includegraphics[width=0.7\columnwidth]{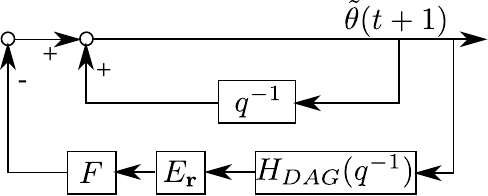}
    \caption{Equivalent feedback representation of the averaged dynamic VS-LMS algorithm.}
    \label{fig_av_eq}
    \end{center}
    \end{figure}

We are now interested on the average of the correcting term $\mathbf{r}(t)\mathbf{r}^T\tilde{\mathbf{w}}(t+1)$. To do averaging, one assumes that the adaptation gain is small enough such that the estimated parameters evolve slowly. This means that:
\begin{multline}
\frac{1}{N+1}\sum_{i=0}^{N}(\mathbf{r}(t-i)\mathbf{r}(t-i)^T\tilde{\mathbf{w}}(t-i+1))  \\
\approx \frac{1}{N+1}\sum_{i=0}^{N}(\mathbf{r}(t-i)\mathbf{r}(t-i)^T)\tilde{\mathbf{w}}(t) 
\end{multline} 
However, for large enough $N$:
\begin{equation}
\frac{1}{N+1}\sum_{i=0}^{N}(\mathbf{r}(t-i)\mathbf{r}(t-i)^T) \approx \textbf{E}\lbrace{\mathbf{r}(t)\mathbf{r}(t)^T}\rbrace=\textbf{E}_{\mathbf{r}}
\end{equation}
and the evolution of the parameter error on the average will be given by:
\begin{align}
\tilde{\mathbf{w}}(t+1)&=\tilde{\mathbf{w}}(t)-FH_{DAG}(q^{-1})[\textbf{E}_{\mathbf{r}}\tilde{\mathbf{w}}(t+1)]\\
&=\tilde{\mathbf{w}}(t)-F\textbf{E}_{\mathbf{r}}H_{DAG}(q^{-1})[\tilde{\mathbf{w}}(t+1)]
\end{align}
leading to the equivalent feedback representation shown in Fig.\ref{fig_av_eq}.
If in addition one makes the assumption that $\mathbf{r}(t)$ is a  stationary persistently exciting signal:
\begin{equation}
\sigma_1 I<\frac{1}{N+1}\sum_{i=0}^{N}(\mathbf{r}(t-i)\mathbf{r}(t-i)^T)<\sigma_2I; ~~ \sigma_1,~\sigma_2>0 
\end{equation}
 it results that $\textbf{E}_{\mathbf{r}}> 0$ (equivalent to the condition of persistent excitation, see \cite{LandauSpringer11}). As a consequence, the equivalent  feedback system will be asymptotically stable (i.e. the parameter error will go to zero as well as the adaptation error) under the sufficient condition that $H_{DAG}(z^{-1})$ is  a strictly positive real (SPR) transfer function  (which is indeed the case) since the feedforward path is passive (an integrator) \cite{LandauSpringer11}.
 Therefore, the stability conditions are relaxed when working with low values of the adaptation gain/step size $\mu>0$.
 
\subsection{Adaptation transient analysis}\label{lin}
\begin{figure}[htb!]%
    \begin{center}
    \includegraphics[width=\columnwidth]{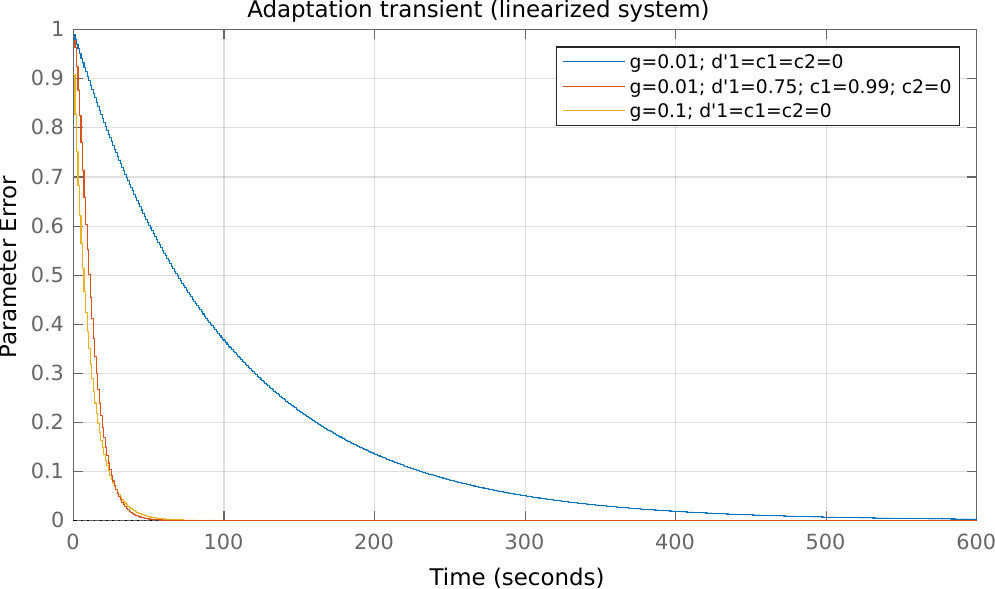}
    \caption{Adaptation transient on the linearized system.}
    \label{fig_lin_eq}
    \end{center}
\end{figure}

 Under the assumptions of persistence of excitation and low adaptation gains (slow adaptation) one can push further the approximation of the equation describing the behavior of the  averaged adaptation algorithm via linearization. This  allows to understand the effect of the DAG and provides further hints for its design.\\
 Consider the case of a single parameter to adapt. Linearization corresponds in Fig. \ref{fig_av_eq} to the replacement of the product $F\textbf{E}_\mathbf{r}$ by a constant positive definite matrix, which for $\dim (\tilde{\mathbf{w}})=1$ is a positive scalar denoted $g$.\\
 The linearized approximation of the algorithm will be described by: 
 \begin{equation}
 \tilde{\mathbf{w}}(t+1)= \tilde{\mathbf{w}}(t)-gH_{DAG}(q^{-1})[\tilde{\mathbf{w}}(t+1)]
 \end{equation}
 which corresponds to a linear feedback system whose output is $\tilde{\mathbf{w}}(t+1)$ (in Fig. \ref{fig_av_eq} one replaces $F\textbf{E}_{\mathbf{r}}$ by $g$). The adaptation transient behavior will be described by the output sensitivity function of this feedback system. For the particular case of an ARIMA2 dynamic adaptation gain, the output sensitivity function will be given by:
 \begin{align}
 S&=\frac{1-q^{-1}}{1+gH_{DAG}(q^{-1})} \nonumber\\
 &=\frac{(1-q^{-1})(1-d'_1q^{-1})}{(1-q^{-1})(1-d'_1q^{-1})+g(1+ c_1q^{-1}+c_2q^{-2})}
\end{align} 
We are interested in the response of this transfer function with respect to a step parameter error.
 Fig.\ref{fig_lin_eq}
shows the step response for various values of the DAG coefficients and two selected values of the adaptation gain $g$.\\
As it can be observed for $g=0.01$ the convergence time for the basic algorithm is approximately 600 s. By adding the DAG  ($d'_1=0.75, c_1=0.99; c2=0$) the convergence time is reduced to approximately 70 s. Note that a similar performance can be obtained with the basic algorithm by multiplying the gain $g$ by 10 (i.e. the acceleration obtained with the DAG without augmenting the adaptation gain is equivalent to using an adaptation gain 10 times larger on the basic algorithm).
\subsection{Stochastic environment -  convergence analysis}
We will consider the I/O model given in \eqref{eq:output} where the output is disturbed by a noise $\mathbf{n}(t+1)$:
\begin{equation}\label{eq:stoch 1}
x(t+1)=\mathbf{w}^T \mathbf{r}(t)+\mathbf{n}(t+1)
\end{equation}
In this equation, $\mathbf{n}(t)$ is a zero mean stationary stochastic disturbance with finite moments.
Using Eqs. \eqref{eq_mu}, \eqref{eq:e-e0} and \eqref{eq:dynamic1}, the PLMS algorithm can be expressed as:
\begin{equation}
\hat{\mathbf{w}}(t+1)= \hat{\mathbf{w}}(t)+ \mu H_{DAG}(q^{-1})[\mathbf{r}(t)e(t+1)]
\label{eq:S2}
\end{equation}
For the analysis of this algorithm in a stochastic context, we will use the ODE (ordinary differential equation) approach of Benveniste-Lam \cite{Benveniste80, Fan88}. Convergence results will be extracted from an ODE which approximate the algorithm by using an average of the correcting term.
We will make the following assumptions: (1) \emph{Stationary processes $\mathbf{r}(t,\hat{\mathbf{w}})$ and $e(t+1,\hat{\mathbf{w}})$ can be defined for $\hat{\mathbf{w}}(t) \equiv \hat{\mathbf{w}}$,} (2) \emph{$\hat{\mathbf{w}}(t)$ generated by the algorithm  belongs infinitely often to the domain $(D_s)$ for which the stationary processes $\mathbf{r}(t,\hat{\mathbf{w}})$ and $\mathbf{n}(t+1),\hat{\mathbf{w}})$ can be defined.}
For the adaptation algorithm given in Eq. \eqref{eq:S2} one has the following result:
\begin{lem}\label{lemstoch} Assume that the Benveniste's smoothness, boundness and mixing conditions \cite{Benveniste80} are satisfied.
\begin{enumerate}
\item Assume that:
\begin{equation}
\sigma_2 I \geq \mathbf{E}\left\{ \mathbf{r}(t, \hat{\mathbf{w}})   \mathbf{r}^T(t , \hat{\mathbf{w}})\right\}\geq \sigma_1 I; ~~\sigma_1,\sigma_2 >0
\label{eq:S4}
\end{equation}
\item Assume that $\mathbf{r}(n,\hat{\mathbf{w}})$ is such (persistence of excitation):
\begin{equation}
\mathbf{r}(n,\hat{\mathbf{w}})^T(\hat{\mathbf{w}}- \mathbf{w}) =0 \Rightarrow \hat{\mathbf{w}}=\mathbf{w}
\label{eq:S5}
\end{equation}
\item Assume that either:
\begin{itemize}
\item $\mathbf{n}(t+1)$ is a sequence of independent equally distributed normal random variables $(0, \sigma)$ or
\item 
\begin{equation}
\mathbf{E}\lbrace \mathbf{r}(t,\hat{\mathbf{w}}),\mathbf{n}(t+1, \hat{\mathbf{w}}\rbrace =0
\end{equation}

\end{itemize}
\item Assume that $H_{DAG}(z^{-1})$ is stricyly positive real
\end{enumerate}
 For a sufficiently small $\mu$ one has: 
\begin{equation}
\mathbf{P}\lbrace \parallel \hat{\mathbf{w}}([T/\mu])-\mathbf{w} \parallel \geq \eta +C\zeta(\mu)\rbrace < C'\zeta(\mu)
\label{eq:S6}
\end{equation}
where: $\zeta(\mu)$ is  a positive decreasing function, $C, C'$ are positive constants, $\overline{\mathbf{w}}(t)$
is the average of $\mathbf{w}(t)$,  $T < \infty $ is such that $\lbrace \parallel \overline{\mathbf{w}}(T)- \mathbf{w} \parallel \rbrace < \eta$ and $[T/\mu]$ is the closest integer to $T/\mu$.
\end{lem}
The proof of this lemma is given in Appendix \ref{proof_stoch1}.\\
The interpretation of this result is as follows: the parametric error (distance) is bounded in \textit{probability }and this bound depends upon the magnitude of $\mu$ \footnote{This is a stronger result than convergence on the average only}.
\subsection{Convergence rate}\label{conv_rate}
 For estimating the convergence rate, one can consider the ODE equation associated with the algorithm and the Lyapunov function candidate used for stability analysis. The ratio $\frac{\vert\dot{V}\vert}{V}$ is an estimation of the convergence rate for large $t$. From Eqs. \eqref{eq:V} and \eqref{eq:Vdot} one gets:
 \begin{equation}
 \frac{\vert\dot{V}\vert}{V}=\frac{(\hat{\mathbf{w}}-\mathbf{w})^T [G_{\mathbf{w}}+G_{\mathbf{w}}^T](\hat{\mathbf{w}}-\mathbf{w})}{(\hat{\mathbf{w}}-\mathbf{w})^T\mu^{-1}(\hat{\mathbf{w}}-\mathbf{w})}
 \label{eq:rate1}
 \end{equation}
 where $ G_{\mathbf{w}}= \mathbf{E}\lbrace H_{DAG}(q^{-1})E_{\mathbf{r}}\rbrace$.
Assuming that the components of the $\mathbf{r} (t)$ are in the low frequency range (this is often the case), then $H_{DAG}(q^{-1})$ can be approximated by its static gain and one has from Eq. \eqref{eq:rate1}:
\begin{equation}
\frac{\vert\dot{V}\vert}{V}\approx \frac{1+\sum_{j=1}^{n_C}c_j}{1-\sum_{j=1}^{n_{D'}}d'_j}\frac{(\hat{\mathbf{w}}-\mathbf{w})^T [\mathbf {E}_{\mathbf{r}}+\mathbf{E}_{\mathbf{r}}^T](\hat{\mathbf{w}}-\mathbf{w})}{(\hat{\mathbf{w}}-\mathbf{w})^T\mu^{-1}(\hat{\mathbf{w}}-\mathbf{w})}
\label{eq:rate}
 \end{equation}
 where $\mathbf {E}_{\mathbf{r}}= \mathbf{E}[\mathbf{r}(t,\hat{\mathbf{w}}), \mathbf{r}(t,\hat{\mathbf{w}})^T]$. Since $H_{DAG}(q^{-1})$ is strictly positive real, its static gain is positive. Provided that the steady state gain $SSG=\frac{1+\sum_{j=1}^{n_C}c_j}{1-\sum_{j=1}^{n_{D'}}d'_j}> 1$ an acceleration with respect to the standard VS-LMS will be obtained and this will be illustrated in the simulation and experimental sections.

\section{Simulation Results}\label{Sim}
\subsection{Adaptive Line Enhancer}
\begin{figure}
\centering
\includegraphics[width=83mm]{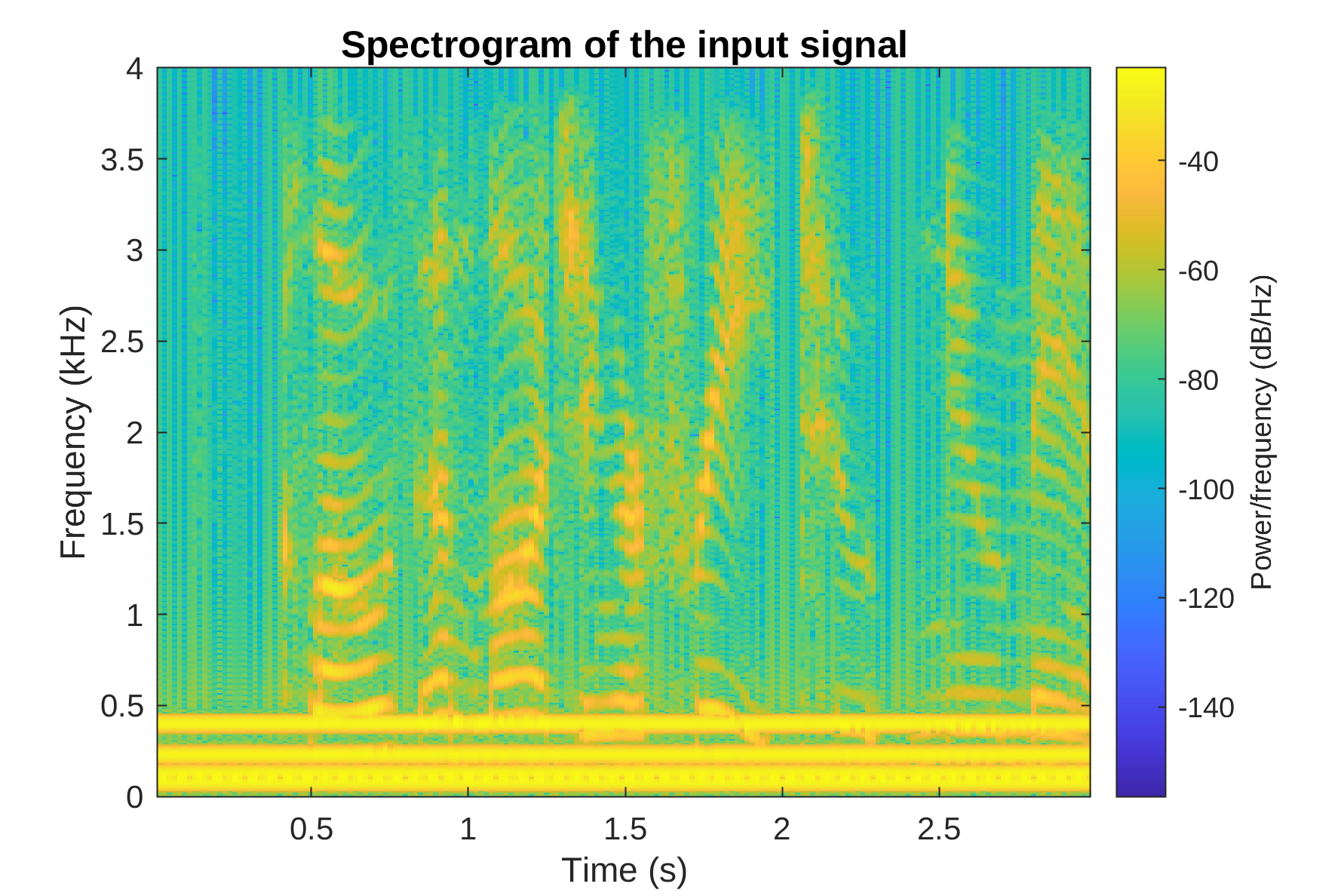}
\caption{The spectrogram of the input signal.}
\label{fig:spectr}
\end{figure}

The properties of several different Dynamic Adaptation Gains (DAGs) were tested in simulations of the Adaptive Line Enhancer (ALE) -- an adaptive filtering
system introduced by Widrow, capable of detecting and removing highly correlated signals embedded in wide-band signals \cite{Widrow75},\cite{Haykin96}. 
The block diagram of the ALE is presented in Fig.~\ref{fig:ale}. The input signal to the ALE, $x(t)$, was constructed of a speech recording with four contaminating 
sine signals. The speech recording was one of the Matlab's Audio Toolbox exemplary files entitled {\em FemaleSpeech-16-8-mono-3secs.wav} (sampling frequency: 8  kHz).
A contaminating signal formed by four sines with frequencies of 80, 125, 230 and 400~Hz has been considered. The signal thus constructed is non-stationary, with the spectrogram presented in Fig.~\ref{fig:spectr}.

\begin{figure}
\centering
\resizebox{6cm}{!}{{\pgfkeys{/pgf/fpu/.try=false}%
\ifx\XFigwidth\undefined\dimen1=0pt\else\dimen1\XFigwidth\fi
\divide\dimen1 by 4139
\ifx\XFigheight\undefined\dimen3=0pt\else\dimen3\XFigheight\fi
\divide\dimen3 by 2147
\ifdim\dimen1=0pt\ifdim\dimen3=0pt\dimen1=8287sp\dimen3\dimen1
  \else\dimen1\dimen3\fi\else\ifdim\dimen3=0pt\dimen3\dimen1\fi\fi
\tikzpicture[x=+\dimen1, y=+\dimen3]
{\ifx\XFigu\undefined\catcode`\@11
\def\temp{\alloc@1\dimen\dimendef\insc@unt}\temp\XFigu\catcode`\@12\fi}
\XFigu8287sp
\ifdim\XFigu<0pt\XFigu-\XFigu\fi
\pgfdeclarearrow{
  name = xfiga2,
  parameters = {
    \the\pgfarrowlinewidth \the\pgfarrowlength \the\pgfarrowwidth\ifpgfarrowopen o\fi},
  defaults = {
	  line width=+7.5\XFigu, length=+120\XFigu, width=+60\XFigu},
  setup code = {
    \dimen7 2.6\pgfarrowlength\pgfmathveclen{\the\dimen7}{\the\pgfarrowwidth}
    \dimen7 2\pgfarrowwidth\pgfmathdivide{\pgfmathresult}{\the\dimen7}
    \dimen7 \pgfmathresult\pgfarrowlinewidth
    \pgfarrowssettipend{+\dimen7}
    \pgfarrowssetbackend{+-1.25\pgfarrowlength}
    \dimen9 -\pgfarrowlength\advance\dimen9 by-0.5\pgfarrowlinewidth
    \pgfarrowssetlineend{+\dimen9}
    \dimen9 -\pgfarrowlength\advance\dimen9 by-0.5\pgfarrowlinewidth
    \pgfarrowssetvisualbackend{+\dimen9}
    \pgfarrowshullpoint{+\dimen7}{+0pt}
    \pgfarrowsupperhullpoint{+-1.25\pgfarrowlength}{+0.5\pgfarrowwidth}
    \pgfarrowssavethe\pgfarrowlinewidth
    \pgfarrowssavethe\pgfarrowlength
    \pgfarrowssavethe\pgfarrowwidth
  },
  drawing code = {\pgfsetdash{}{+0pt}
    \ifdim\pgfarrowlinewidth=\pgflinewidth\else\pgfsetlinewidth{+\pgfarrowlinewidth}\fi
    \pgfpathmoveto{\pgfqpoint{-1.25\pgfarrowlength}{-0.5\pgfarrowwidth}}
    \pgfpathlineto{\pgfqpoint{0pt}{0pt}}
    \pgfpathlineto{\pgfqpoint{-1.25\pgfarrowlength}{0.5\pgfarrowwidth}}
    \pgfpathlineto{\pgfqpoint{-\pgfarrowlength}{0pt}}
    \pgfpathclose
    \ifpgfarrowopen\pgfusepathqstroke\else\pgfsetfillcolor{.}
	\ifdim\pgfarrowlinewidth>0pt\pgfusepathqfillstroke\else\pgfusepathqfill\fi\fi
  }
}
\clip(-22,-2307) rectangle (4117,-160);
\tikzset{inner sep=+0pt, outer sep=+0pt}
\pgfsetarrows{[line width=7.5\XFigu]}
\pgfsetarrowsend{xfiga2}
\pgfsetlinewidth{+7.5\XFigu}
\pgfsetdash{{+60\XFigu}{+60\XFigu}}{++0pt}
\pgfsetstrokecolor{black}
\draw (1842,-1798)--(1842,-1663)--(2517,-673);
\pgfsetdash{}{+0pt}
\draw  (3375,-219) circle [radius=+51];
\pgfsetlinewidth{+15\XFigu}
\draw (0,-225)--(3330,-225);
\draw (3420,-225)--(4095,-225);
\draw (675,-225)--(675,-540);
\pgfsetlinewidth{+7.5\XFigu}
\pgfsetarrowsend{}
\draw (675,-1035)--(675,-1440);
\draw (405,-540) rectangle (945,-1035);
\pgfsetlinewidth{+15\XFigu}
\pgfsetarrowsend{xfiga2}
\draw (675,-1260)--(1530,-1260);
\draw (2610,-1260)--(3375,-1260)--(3375,-270);
\pgfsetlinewidth{+7.5\XFigu}
\pgfsetfillcolor{white}
\filldraw (1530,-1800) rectangle (2610,-2295);
\draw (675,-1260)--(675,-2070)--(1530,-2070);
\draw (3645,-225)--(3645,-2070)--(2610,-2070);
\filldraw (1530,-990) rectangle (2610,-1485);
\pgfsetfillcolor{black}
\pgftext[base,left,at=\pgfqpointxy{90}{-450}] {\fontsize{24}{26.4}$x(t)$}
\pgftext[base,left,at=\pgfqpointxy{3465}{-405}] {\fontsize{24}{33.6}$-$}
\pgftext[base,left,at=\pgfqpointxy{3150}{-405}] {\fontsize{24}{28.8}\normalfont +}
\pgftext[base,left,at=\pgfqpointxy{990}{-1665}] {\fontsize{24}{26.4} $d(t)$}
\pgftext[base,left,at=\pgfqpointxy{585}{-855}] {\fontsize{24}{26.4}$q^{-\Delta}$}
\pgftext[base,left,at=\pgfqpointxy{1710}{-2025}] {\fontsize{24}{24}\normalfont\rmfamily Adaptation}
\pgftext[base,left,at=\pgfqpointxy{1755}{-2205}] {\fontsize{24}{24}\normalfont\rmfamily algorithm}
\pgftext[base,left,at=\pgfqpointxy{2745}{-1125}] {\fontsize{24}{24}$\hat{y}^\circ(t)$}
\pgftext[base,left,at=\pgfqpointxy{3690}{-450}] {\fontsize{24}{26.4} $ e^\circ (t)$}
\pgftext[base,left,at=\pgfqpointxy{1710}{-1305}] {\fontsize{24}{28.8} $\hat{W}(q^{-1})$}
\endtikzpicture}%
}
\caption{Block diagram of the Adaptive Line Enhancer.}
\label{fig:ale}
\end{figure}
The adaptive filter length and the decorrelation delay ($\Delta$) used in the simulations were both 100. The FIR filter was longer than can be expected, provided that only four
sine signals are to be detected. However, it is known that when the sines are close to each other, the required filter length to detect and cancel them efficiently is 
usually higher.\footnote{This effect was confirmed in simulations with shorter filters, which were not capable of cancelling all the four sines}. The adaptive filter was always started with zero initial conditions.

Both the NLMS and PLMS algorithms, incorporating a DAG, were tested and compared in the simulations, while the NLMS and PLMS without the DAG  served as a reference for comparison.
The step size for the NLMS was experimentally selected as 0.02. The step size for the PLMS ($5\cdot 10^{-4}$) was adjusted to achieve the same initial convergence speed for both 
the NLMS and PLMS algorithms without the DAG.

\begin{figure}
\centering
\includegraphics[width=83mm]{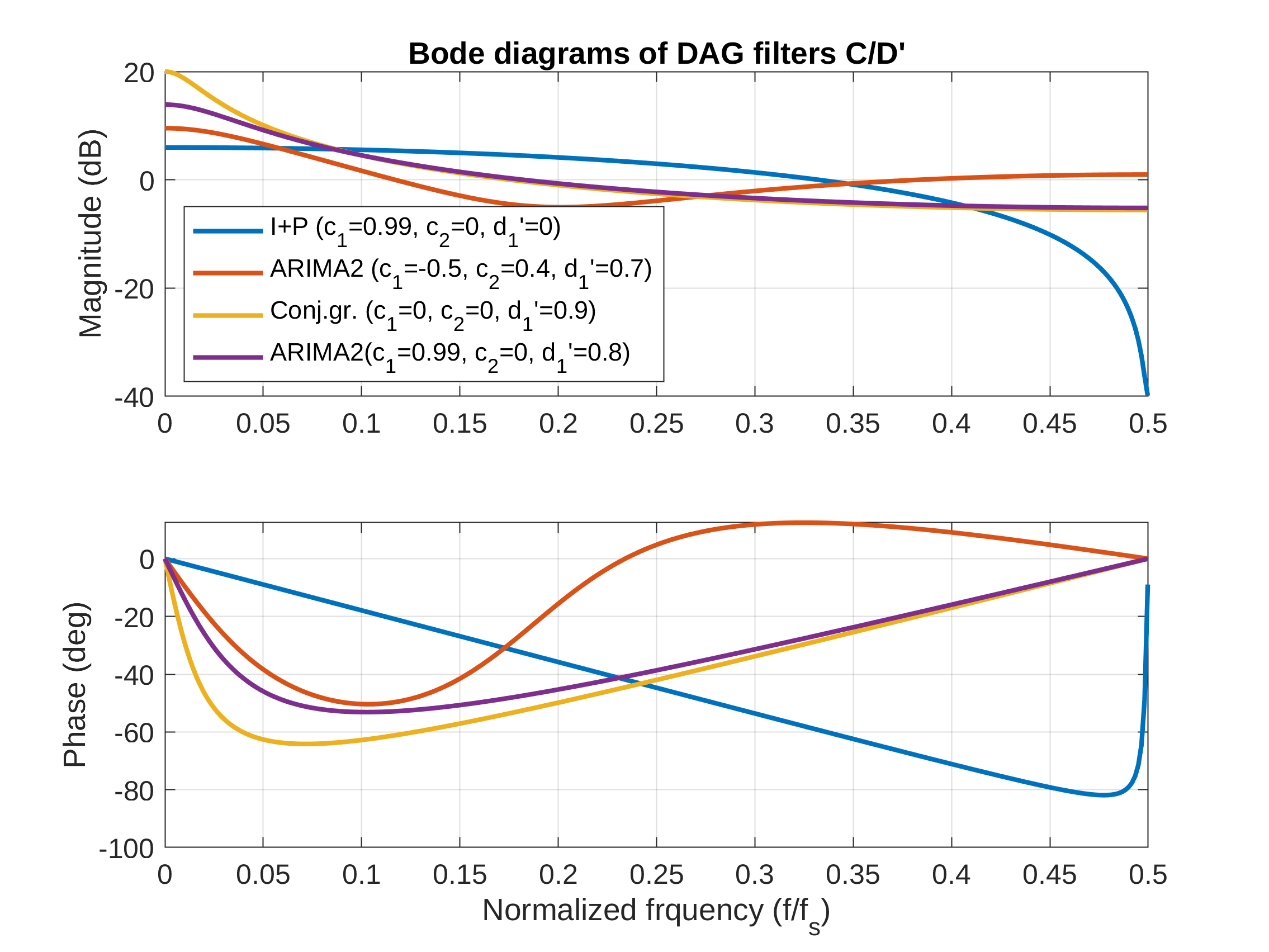}
\caption{Bode diagrams of the DAGs used in simulations of the ALE.}
\label{fig:dagbode}
\end{figure}

The performance of different DAG settings were evaluated by the means of two quantities. The first was the convergence speed calculated as time (in samples)
after which the mean squared error (MSE) achieved the level of $-40$~dB for the first time. The second was the sum of the MSE during the first 3200 iterations (0.4 s) of the simulation.
 Please note that excess MSE could not have been used due to the lack of the Wiener filter model
for this nonstationary input signal, which implies impossibility to calculate the system noise. Fig. \ref{fig:dagbode} summarizes the frequency characteristics of the various DAG used. All are SPR.

\begin{figure}
\centering
\includegraphics[width=83mm]{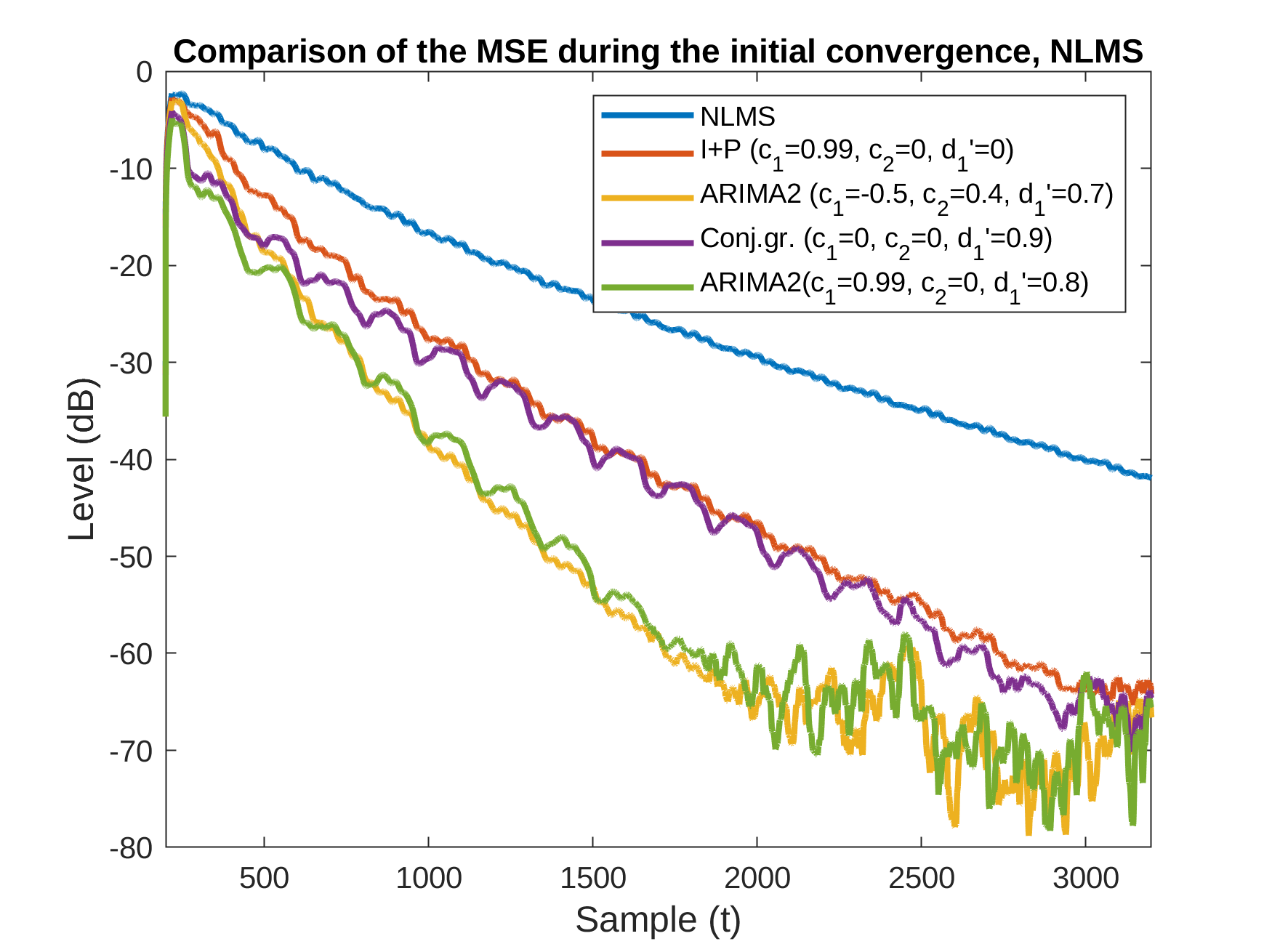}
\caption{The spectrogram of the input signal.}
\label{fig:msenlms}
\end{figure}

Figure \ref{fig:msenlms} presents the results obtained with the NLMS algorithm for the selected set of coefficients of the $H_{DAG}(q^{-1})$. The figure is zoomed to show the MSE during
the initial convergence of the ALE filter. The plots clearly show that all the selected  DAGs speed up the initial convergence significantly, compared to the NLMS without the DAG.
This is confirmed by numerical values of the selected performance indicators, presented in Table~\ref{tab:nlms}.
In case of the NLMS algorithm, the fastest convergence has been obtained for  set no. 5 ($c_1=-0.5$, $c_2=0.4$ and $d_1'=0.7$). In the same time, the smallest sum of MSE is obtaned for set no. 6
($c_1=0.99$, $c_2=0$, $d_1'=0.8$). 

\begin{figure}
\centering
\includegraphics[width=83mm]{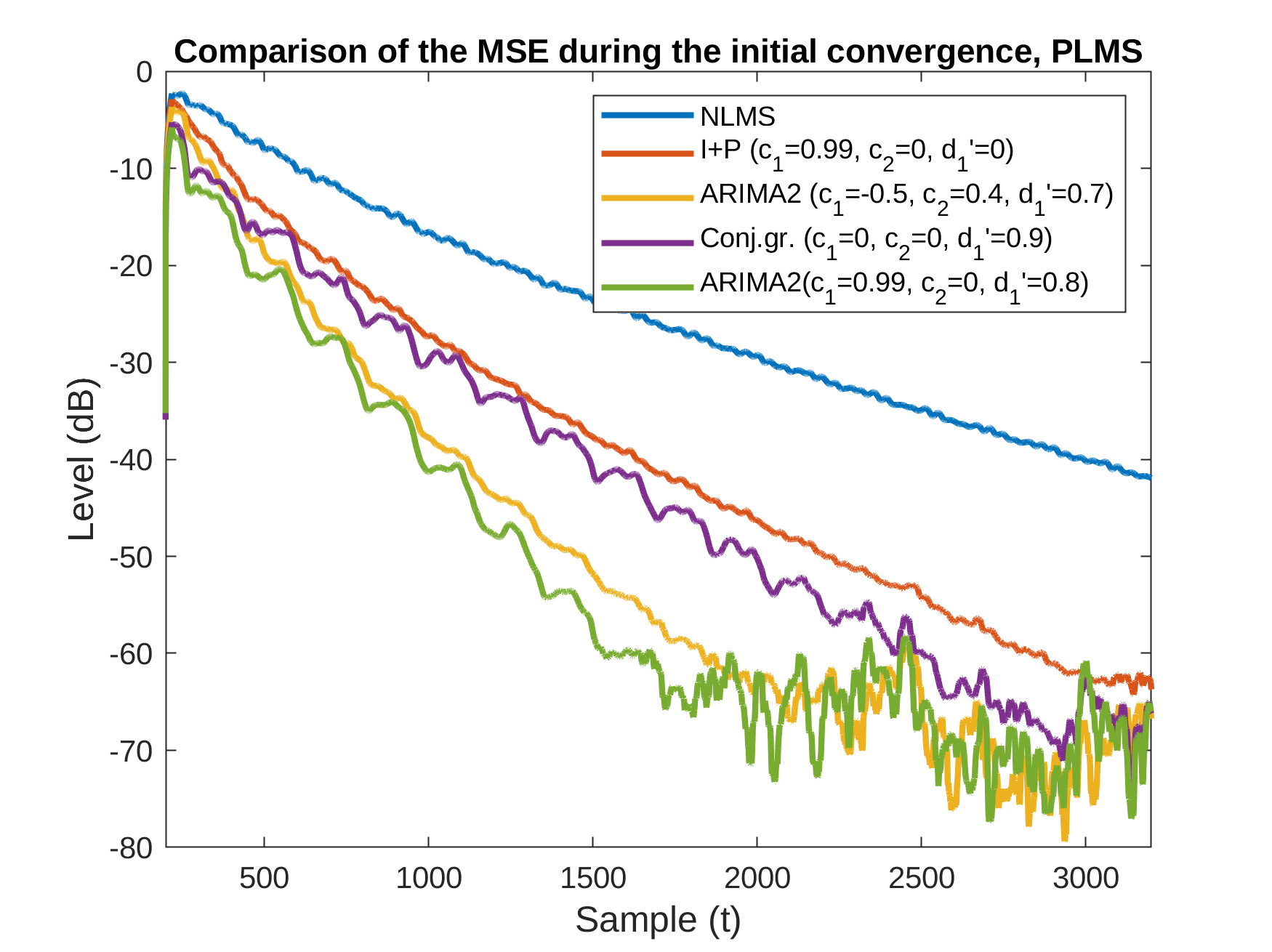}
\caption{The spectrogram of the input signal.}
\label{fig:mseplms}
\end{figure}

Figure~\ref{fig:mseplms} presents similar results obtained for the PLMS algorithm, and the corresponding values of the performance indicators are presented in Table~\ref{tab:plms}.
In the case of this algorithm, set no. 6 provides both the shortest convergence time and the smallest sum of MSE.

\begin{table}
\caption{Results of the ALE experiments averaged over 250 runs, ALE length: 100}
\label{tab:nlms}

\footnotesize 
\begin{tabular}{||l|l|c|c||}
\hline
\multirow{2}{*}{No} & \multirow{2}{*}{Algorithm / DAG settings} & Conv.time & Sum of \\ 
& & (samples) & MSE \\ \hline
1 & NLMS & 2791 & 145.2 \\ \hline
3 & $c_1=0.99, c_2=0, d_1'=0$ & 1426 & 74.4 \\ \hline
4 & $c_1=0, c_2=0, d_1'=0.9$ & 1296 & 33.7 \\ \hline
5 & $c_1=-0.5, c_2=0.4, d_1'=0.7$ & 843 & 49.8 \\ \hline
6 & $c_1=0.99, c_2=0, d_1'=0.8$ & 916 & 24.9 \\ \hline
\end{tabular}

\end{table}

\begin{table}
\caption{Results of the ALE experiments averaged over 250 runs, PLMS algorithms.}
\label{tab:plms}

\footnotesize 
\begin{tabular}{||l|l|c|c||}
\hline
\multirow{2}{*}{No} & \multirow{2}{*}{Algorithm / DAG settings} & Conv.time & Sum of \\ 
& & (samples) & MSE \\ \hline
1 & NLMS & 2791 & 145.2 \\ \hline
2 & PLMS & 2796 & 132.0 \\ \hline
3 & $c_1=0.99, c_2=0, d_1'=0$ & 1429 & 62.1 \\ \hline
4 & $c_1=0, c_2=0, d_1'=0.9$ & 1286 & 31.5 \\ \hline
5 & $c_1=-0.5, c_2=0.4, d_1'=0.7$ & 914 & 41.5 \\ \hline
6 & $c_1=0.99, c_2=0, d_1'=0.8$ & 770 & 20.5 \\ \hline
\end{tabular}

\end{table}

\subsection{Filter identification}

Assume that the signal x(t) in Fig. \ref{fig_LMS_basic_scheme} is the output of an IIR filter whose input is a PRBS (pseudo random binary sequence) with $N=2^8-1=255$ samples. The IIR filter is characterized by the transfer operator (unknown).
\begin{equation}
Sys=\frac{q^{-2}+0.5q^{-3}}{1-1.5q^{-1}+0.7q^{-2}}
\label{eq:Syssim}
\end{equation}
 The objective was to estimate the parameters of this IIR filter as well as the parameters of an FIR filter able to approximate this IIR filter. The PLMS algorithm has been used to illustrate the properties of the various DAGs. 
\subsubsection{IIR identification}\label{pf}
For the  adaptation gain/learning rate $\mu=0.02$, the performance of the adaptation algorithms was evaluated with respect to the choice of the coefficients $c_1, c_2, d'_1$. To assess the performance, the following indicators were used: (i) the sum of the squared \emph{a posteriori} prediction errors:
$J_{\epsilon}(t)=\sum_{t=0}^{N}e^2(t+1)$, 
(ii) the square of the parametric distance:
$D^{2}(t)=\left\{[\mathbf{w}-\hat{\mathbf{w}}(t)]^T[\mathbf{w}-\hat{\mathbf{w}}(t)]\right\}$, 
and (iii) the sum of the  squared parametric distance:
$J_{D}(t)=\sum_{t=0}^{N}D^{2}(t)$.
\begin{figure}[htb]
    \includegraphics[width=1\columnwidth]{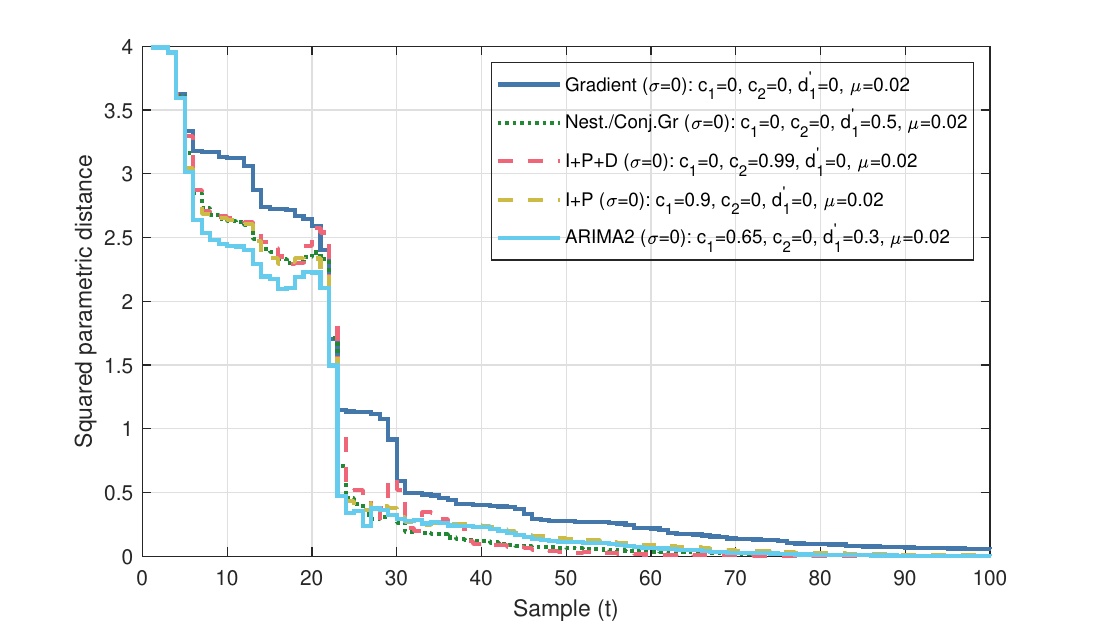}
		\includegraphics[width=1\columnwidth]{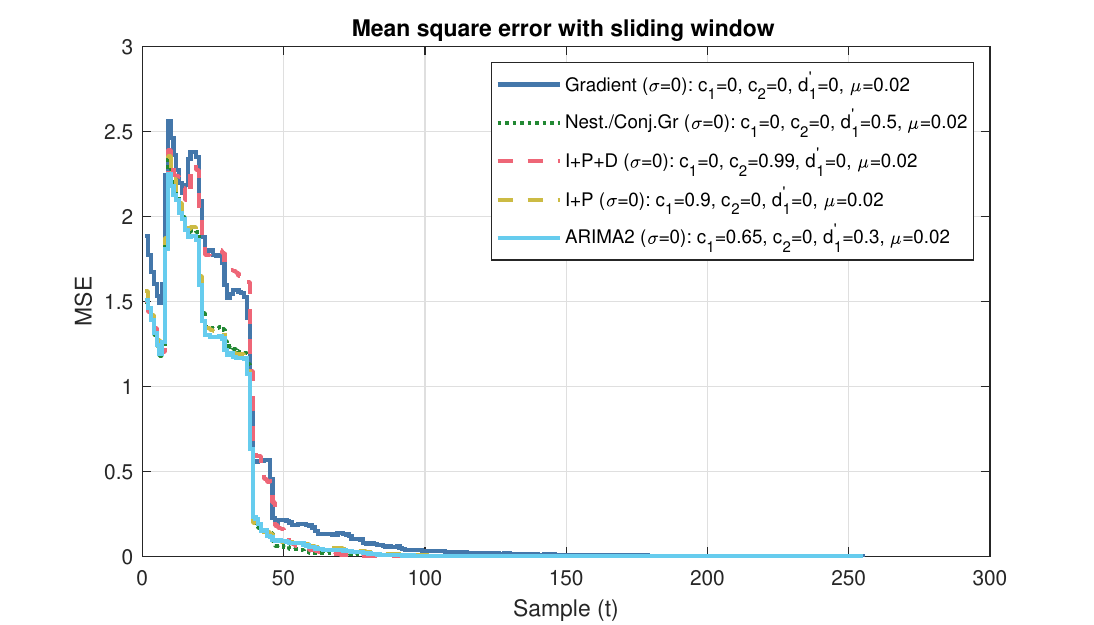}
	\caption{Evolution of the squared parametric distance $D^2(t)$ (top) and of the MSE (bottom)(IIR identification).}
	\label{fig:Distance}
\end{figure}
Table \ref{perf} summarizes the performance of the 2nd order ARIMA algorithm (ARIMA 2) and of the various particular cases. All the DAGs in  Table \ref{perf} are SPR.
  Fig.~\ref{fig:Distance} shows the evolution of the squared  parametric distance (top) and of the MSE (mean square error)(bottom). Clearly the use of a DAG  provided a significant performance improvement with respect to the basic PLMS algorithm. The ARIMA 2 DAG gave the best results.

 \renewcommand{\tabcolsep}{1pt}
  \begin{table}[htb!]
    \begin{center}%
    \caption{Performance of 2nd order ARIMA algorithms.}
        \begin{tabular}{|c|c|c|c|c|c|c|c|}\hline
            Algorithm   &$c_1$  		& $c_2$ & $d'_1$ & $J_{D}(N)$ &	$J_{\epsilon}(N)$&$J_{D}(N)$&$J_{\epsilon}(N)$
            \\ \hline
               &  		&  &  & IIR &	IIR & FIR & FIR  \\ \hline \hline
            Gradient  &$0$ 	& $0$ 	& $0$ 	& $\textbf{91.23}$ & $\textbf{81.44}$&\textbf{589.1}&\textbf{396.7}
            \\ \hline
            Conj.Gr/Nest..& $0$ 	& $0$ 	& $0.5$ 	& $66.95$ & $60.27$&350.3&249.6  
            \\ \hline
            I+P+D  & $0$ 	& $0.99$ 	& $0$ 	& $68.86$ & $74.15$ &351.7&255
            \\ \hline
            I+P & $0.9$ &  $0$ & $0$ & $69.71$ & $61.66$& 356.8&247.6
            \\ \hline
            ARIMA 2 & $0.65$ &  $0$ & $0.3$ & $\textbf{64.45}$ & $\textbf{59.93}$&\textbf{ 340.9}&\textbf{222}
            \\ \hline
        \end{tabular}
        \label{perf}
    \end{center}  
\end{table}

\begin{figure}[ht!]
		\centering
		\includegraphics[width=\columnwidth]{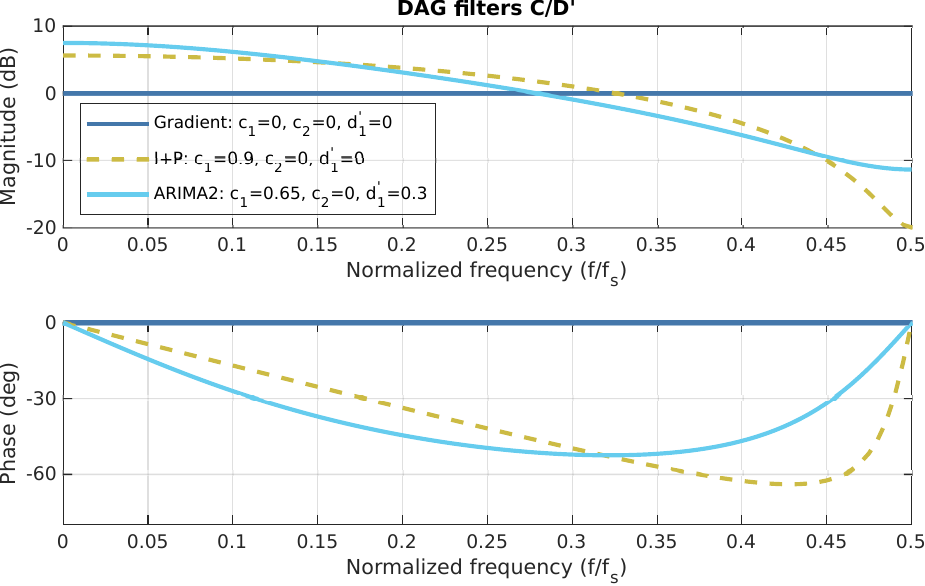}
		\caption{Bode diagram of the dynamic adaptation gain/learning rate $H_{DAG}$ for ARIMA2, I+P and Gradient  algorithms (from Table \ref{perf1}).}
		\label{fig:C_D'}
	\end{figure}

 Fig. \ref{fig:C_D'} gives the Bode diagram for the ARIMA 2 and I+P algorithms (the gradient algorithm corresponds to the $0$~dB axis). One can notice that the phase lag is less than 90 degrees at all the frequencies. It was verified that the average gain over the all frequency range was $0$~dB. This means that the improvement in performance is related to the frequency distribution of the adaptation gain/learning rate. Now examining the magnitude, one observes that both are low pass filters amplifying low frequencies. This means that if the frequency content of the gradient is in the low frequency range, the effective adaptation gain/learning rate will be larger than the gradient adaptation gain (0 dB), which should have a positive effect upon the adaptation/learning transient. In particular, the DAG which has a larger gain in low frequencies (ARIMA2) should provide better performance than the (I+P) DAG (which was indeed the case). This observation is also coherent with the estimated asymptotic convergence rate.\\
\subsubsection{FIR identification}
\begin{figure}[htb]
    \includegraphics[width=1\columnwidth]{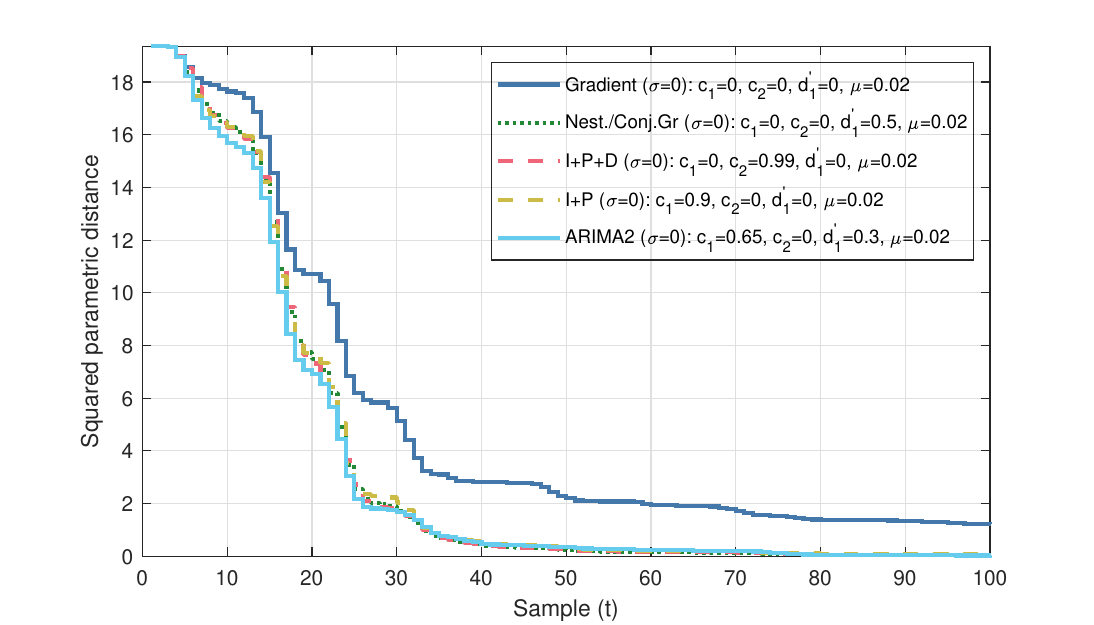}
		\includegraphics[width=1\columnwidth]{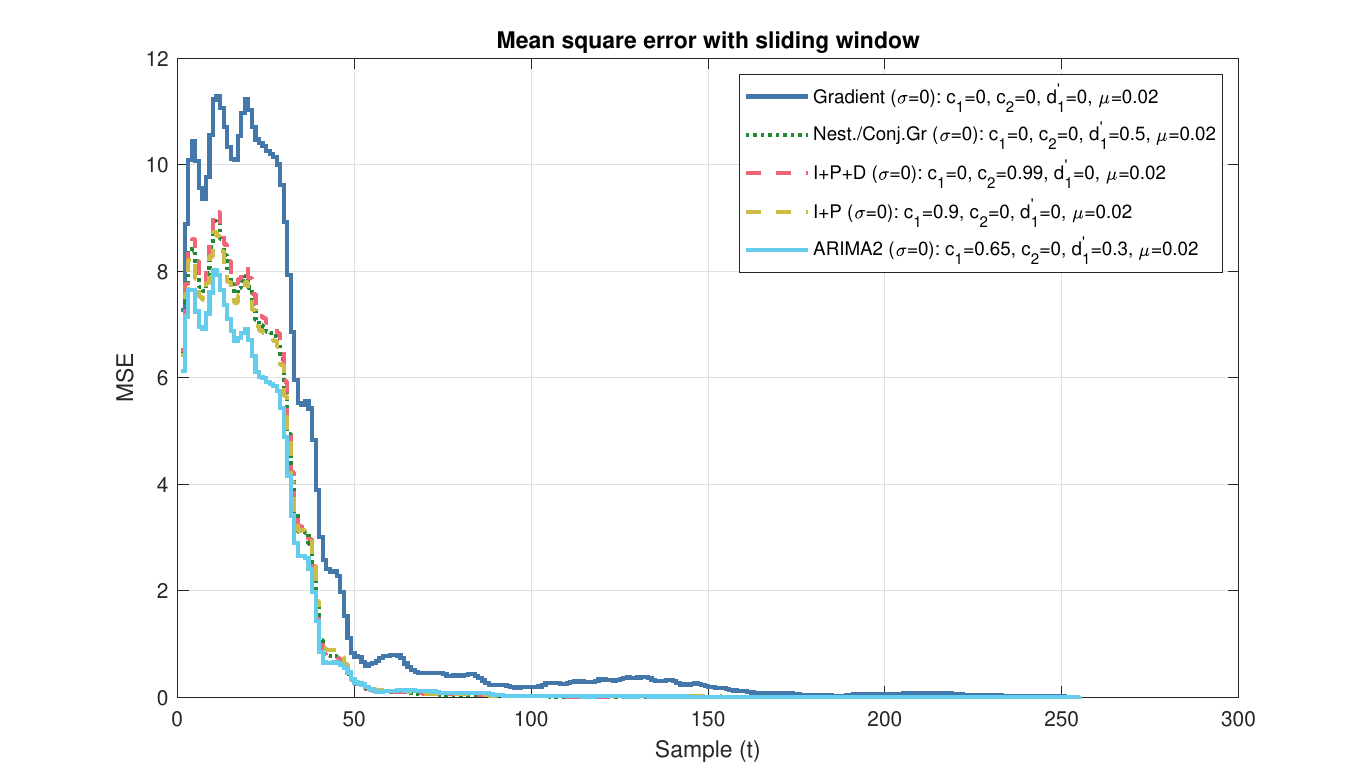}
	\caption{Evolution of the squared parametric distance $D^2(t)$ (top) and of the MSE (bottom) for the identification of a FIR with 30 parameters.}
	\label{fig:impmatching}
\end{figure}

A FIR filter with 30 parameters allows to well approximate the IIR filter given in Eq. \eqref{eq:Syssim}. To evaluate the adaptation transient for the estimated parameters an estimation of "good" parameters was done first.
Figure \ref{fig:impmatching} gives  the time evolution of the squared parametric distance (top) and of the MSE (bottom). Performance of the various DAGs are summarized in Table~\ref{perf} (last two columns).
  The conclusions drawn for the IIR filter identification hold also for FIR filter identification.
\subsubsection{Stochastic case}
\begin{figure}[ht!]
		\centering
		\includegraphics[width=\columnwidth]{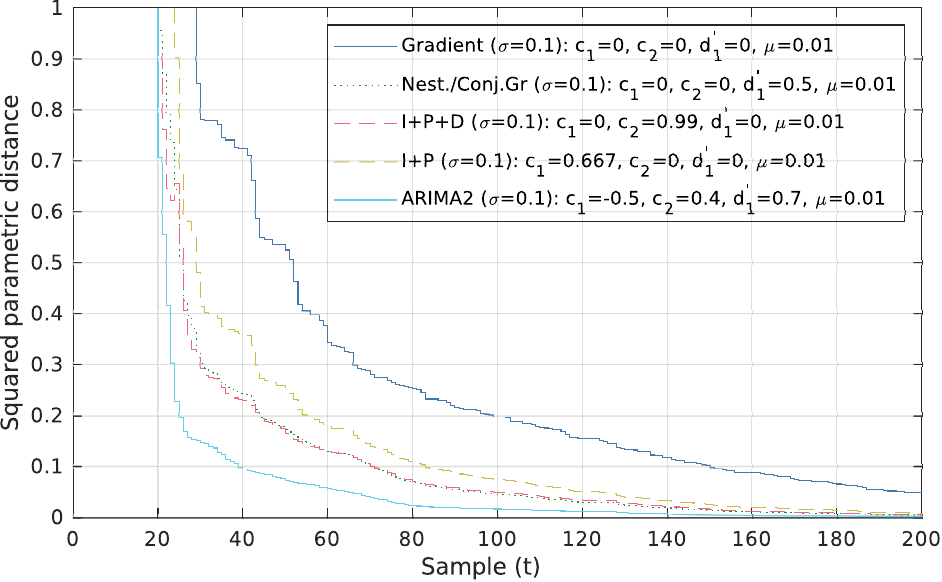}
		\caption{Evolution of the squared parametric distance in the presence of noise (zoom).}
		\label{fig:stoch1}
	\end{figure}

To the same simulation example, a white noise was added on the output (signal/noise ratio (standard deviation): 33 dB). The adaptation gain/step size  $\mu=0.01$ was used. The input in this case was a PRBS with $N=2^{11}-1=20427$ samples. Figure \ref{fig:stoch1} shows a zoom of the evolution of the squared parametric distance (average over $100$ noise realizations) for IIR filter identification. One gets almost asymptotic unbiased parameters estimates (initial value of the squared parametric distance is 4) and the improvement of the transient performances with respect to the basic algorithm is obvious\footnote{The transient performance can be related to the asymptotic convergence rate  given in Section \ref{conv_rate}.}.

\section{Experimental results}\label{Exp}
\begin{figure}[htb!]
  \centering
  
		\includegraphics[width=\columnwidth]{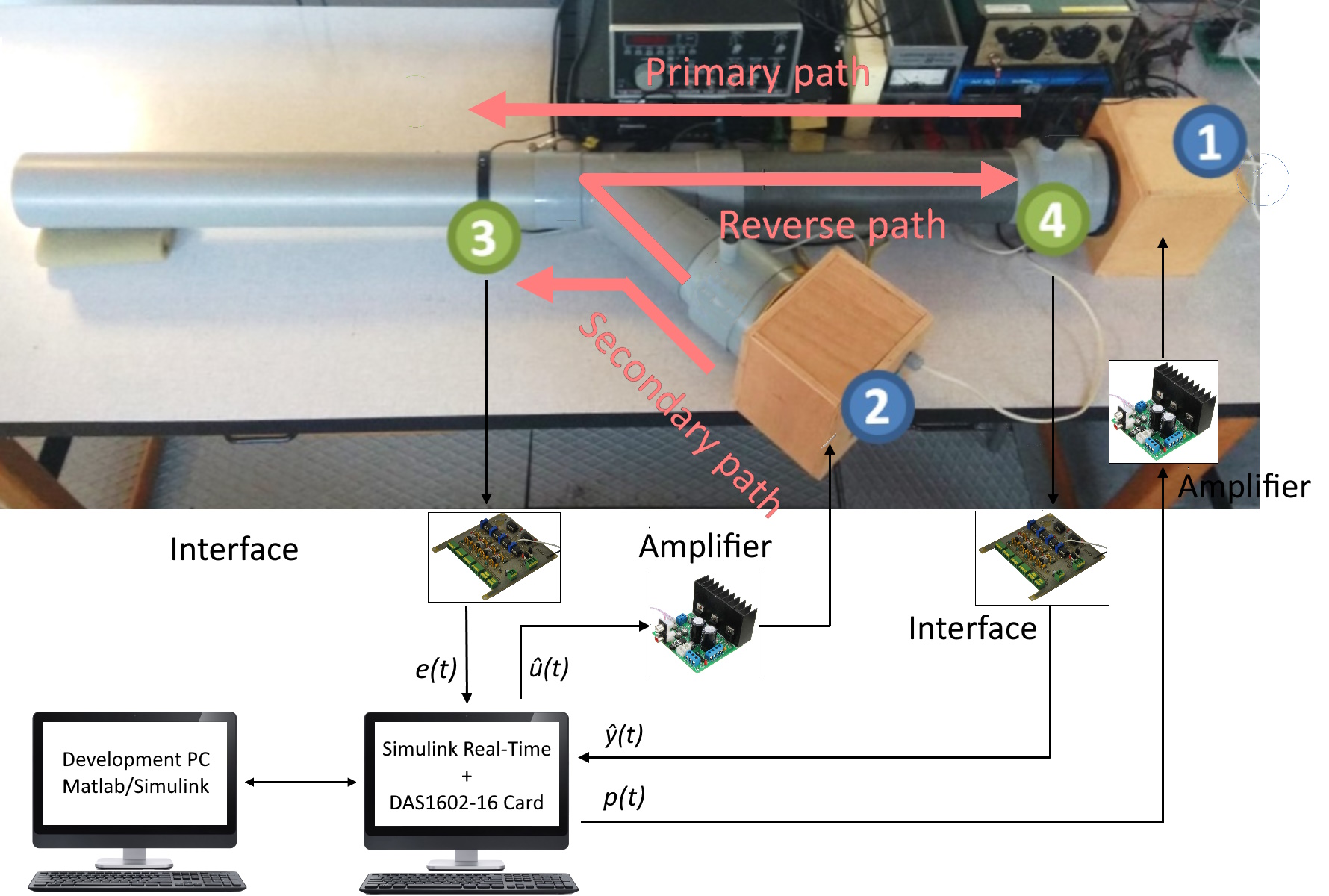} %
	\caption{Duct active noise control test-bench.}
	\label{photos_test_bench}
\end{figure}

The objective of his section is to show that the dynamic adaptation gain can be implemented on various VS-LMS algorithms and this will lead to a significant acceleration of the adaptation transient. Specifically, for this paper the LMS, NLMS and PLMS algorithms were implemented and tested experimentally on an active noise control test-bench (adaptive feedforward noise attenuation).
 Figure~\ref{photos_test_bench} shows the view of the test-bench used for experiments. Detailed description can be found in \cite{LandauJSV23}.
The speakers and the microphones were connected to a target computer with Simulink Real-Time\textsuperscript{\textregistered}. A second computer is used for development and operation with Matlab/Simulink. The sampling frequency was $f_s=2500$~Hz. 

The various paths are described by models of the form:
$X(q^{-1})=q^{-d_x}\frac{B_{X}(q^{-1})}{A_{X}(q^{-1})} =q^{-d_x}\frac{b^{X}_{1}q^{-1}+...+b^{X}_{n_{B_{X}}}q^{-n_{B_{X}}}}{1+a^{X}_{1}q^{-1}+...+a^{X}_{n_{A_{X}}}q^{-n_{A_{X}}}}$, 
with $B_{X}=q^{-1}B^*_{X}$ for any $X \in \{G,M,D\}$. 
$\hat{G}=q^{-d_G}\frac{\hat{B}_G}{\hat{A}_G}$, $\hat{M}=q^{-d_M}\frac{\hat{B}_M}{\hat{A}_M}$, and
$\hat{D}=q^{-d_D}\frac{\hat{B}_D}{\hat{A}_D}$ denote the identified (estimated) models of the secondary ($G$), reverse ($M$), and primary ($D$) paths. The system's order is defined as (the indexes $G$, $M$, and $D$ have been omitted): $n= max(n_A, n_{B}+d)$.

The models of the above paths are characterized by the presence of many pairs of very low damped poles and zeros. These models were identified experimentally. The  orders of the various identified models were: $n_G=33$, $n_M=27$ and $n_D=27$.

\begin{figure}[ht]%
    \begin{center}
    \includegraphics[width=8cm]{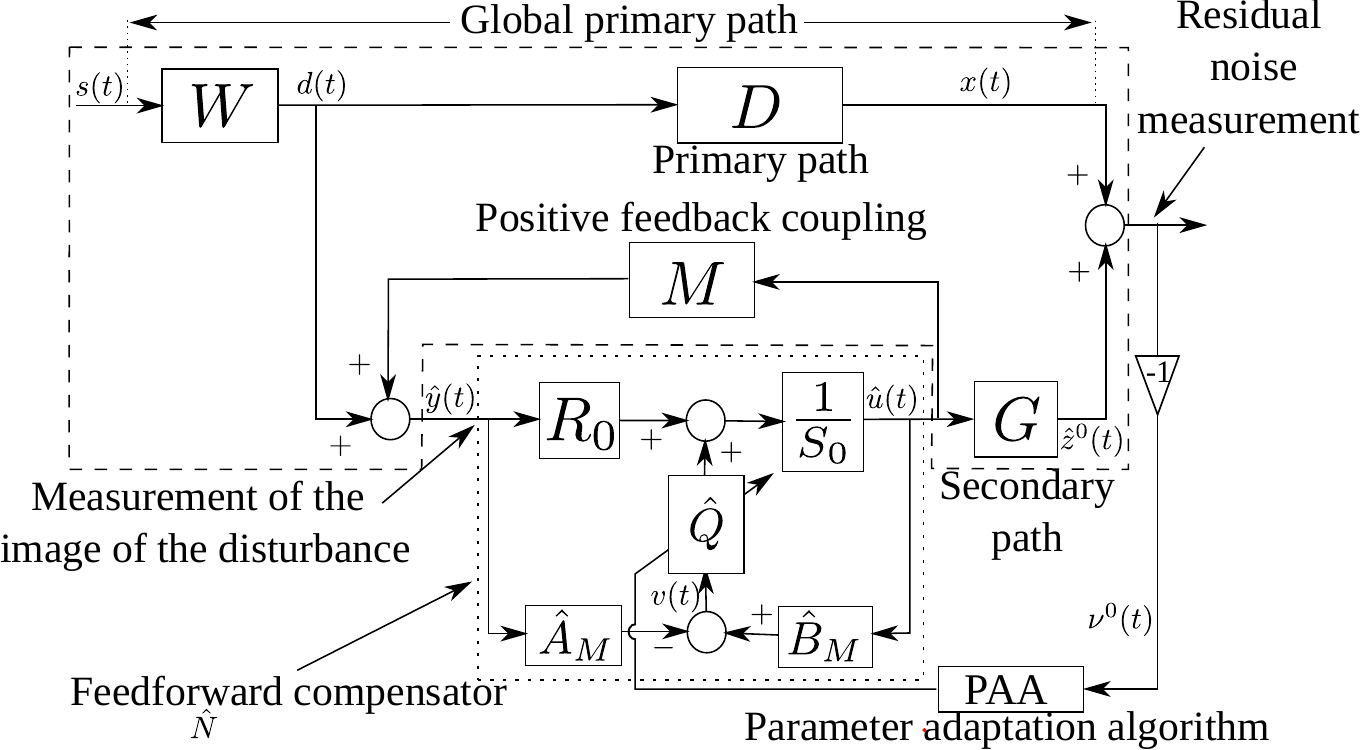}
    \caption{Feedforward ANC with FIR-YK adaptive feedforward compensator.}
    \label{fig_feedforwardAVC}
    \end{center}
\end{figure}

The objective is to attenuate an incoming unknown broad-band noise disturbance eventually mixed with several tonal disturbances. The corresponding block diagram for the adaptive feedforward noise compensation using FIR Youla-Kucera (FIR-YK) parametrization of the feedforward compensator (introduced in \cite{LandauAuto12} for active vibration control and in \cite{airimitoaie:hal-02947816} for active noise control) is shown in Figure~\ref{fig_feedforwardAVC}.

The adjustable filter $\hat{Q}$ has the FIR structure:
\begin{equation}\label{eq_Qhat(q-1)}
            \hat{Q}(q^{-1}) = {\hat{q}_{0}}+\hat{q}_{1}q^{-1}+...+\hat{q}_{n_{Q}}q^{-n_{Q}}
        \end{equation}
and the parameters $\hat{q}_{i}$ will be adapted in order to minimize the residual noise.

The algorithm which was used (introduced in \cite{LandauAuto12}) can be summarized as follows.
One defines:
\begin{eqnarray}
\mathbf{w}^T&=&[q_0,q_1,q_2,\dots,q_{n_Q}]\label{teta}\\
		 \hat{\mathbf{w}}^T&=&[\hat{q}_0,\hat{q}_1,\hat{q}_2,\dots,\hat{q}_{n_Q}]\label{tetahat}\\
		 \mathbf{v}^T(t)&=&[v(t+1),v(t),\dots,v(t-{n_Q}+1)]\label{fi}
	\end{eqnarray}
where:
\begin{equation}
		 v(t+1)=B_{M}\hat{u}(t+1)-A_{M}\hat{y}(t+1)=B^{*}_{M}\hat{u}(t)-A_{M}\hat{y}(t+1)\label{eq:alpha}
	\end{equation}
	One defines also the regressor vector (a filtered observation vector) as:
 \begin{equation}\label{eq_phif}
\mathbf{r}(t)=L(q^{-1})\mathbf{v}(t)=[v_f(t+1),v_f(t),\dots,v_f(t-n_Q+1)]
\end{equation}
where
\begin{equation}
v_{f}(t+1)=L(q^{-1})v(t+1)\label{eq:filtalg2}
\end{equation}
Using $R_0=0$ and $S_0=1$, the poles of the internal positive closed loop are defined by $A_M$ and they will remain unchanged. The filter  used in \eqref{eq:filtalg2} becomes $L=\hat{G}$ and the associated linear transfer operator appearing in the feedforward path of the equivalent feedback system is:
\begin{equation}\label{eq_Halgo2}
H(q^{-1})=\frac{G(q^{-1})}{\hat{G}(q^{-1})}
\end{equation}
(the algorithm uses an approximate gradient). The transfer function associated to $H(q^{-1})$ should be SPR in order to assure asymptotic stability in the case of perfect matching \cite{LandauAuto12}. This is a very mild condition as long as a good experimental identification of the models is done.

The VS-LMS algorithms which were used are of the form given in \eqref{eq_vslms}, where $e^\circ(t)$ is the measured residual noise with minus sign ($\nu^\circ (t)$ in Fig. \ref{fig_feedforwardAVC}),
$\hat{\mathbf{w}}$ is given by Eq. \eqref{tetahat} and $\mathbf{r}$ is given in  Eq. \eqref{eq_phif}\footnote{In signal processing literature when using a filtering of the observation vector, the algorithms are termed FU-VSLMS.}.  The adjustable filter $\hat{Q}$ had 60 parameters.\\
The values of the adaptation gain $\mu$ for the three algorithms were tuned such that in the absence of the DAG, performance was close for the three algorithms. A specificity of this application is also the low value of the average of $\mathbf{r}(t)^T\mathbf{r}(t) ~~(<<1)$. This means that the LMS and PLMS show a very close behavior for a given $\mu$ and that the adaptation gain $\mu$ for the NLMS should be much lower in order to get similar performances. The choices of $\mu$ were 0.2 for LMS, 0.0003 for NLMS and 0.22 for the PLMS.

\renewcommand{\tabcolsep}{4pt}
 \begin{table}[htb!]
    \begin{center}%
    \caption{Parameters of ARIMA2 Dynamic Adaptation Gain.}
        \begin{tabular}{c|c|c|c|c|c}\hline
            Algorithm &$H_{PAA}$--PR &$H_{DAG}$--SPR&$c_1$  		& $c_2$ & $d'_1$ 
            \\ \hline \hline
            Integral (gradient) & Y &Y&$0$ 	& $0$ 	& $0$ 	
            \\ \hline
            Conj.Gr/Nest..& N&Y& $0$ 	& $0$ 	& $0.9$ 	
            \\ \hline
            I+P+D & N &Y& $1.4$ &  $0.5$ & $0$  
            \\ \hline
            I+P & Y&Y & $0.99$ &  $0$ & $0$ 
            \\ \hline
            ARIMA 2 & N&Y& $0.99$ &  $0$ & $0.9$ 
            \\ \hline
        \end{tabular}
        \label{perf1}
    \end{center}  
\end{table}
\begin{figure}[htb]
  \includegraphics[width=1\columnwidth]{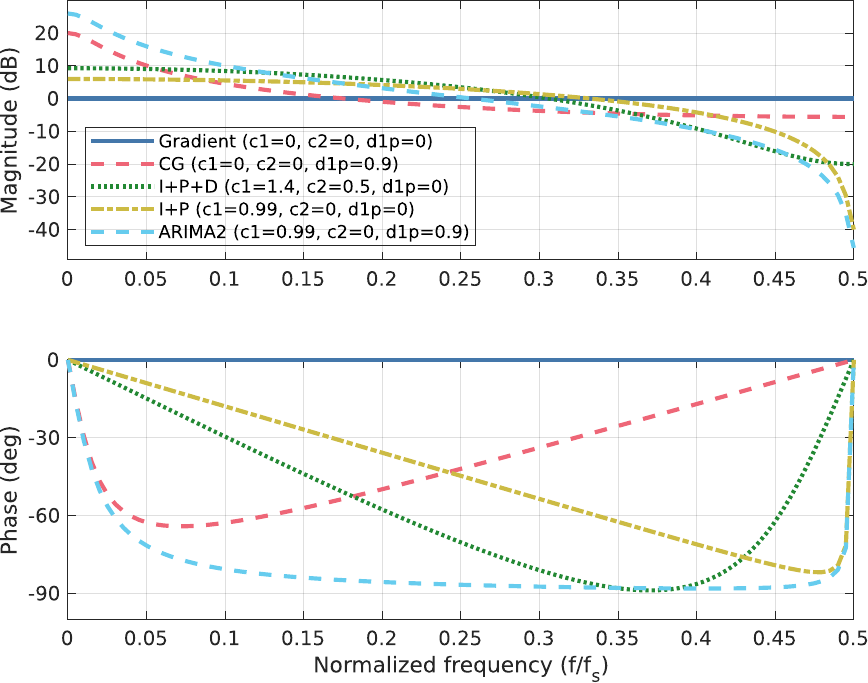}
  \caption{Frequency characteristics of DAGs used in the experiments (see also Table~\ref{perf1}).}
  \label{fig_bode_comp}
\end{figure}

\begin{figure}[htb!]   
    \includegraphics[width=1\columnwidth]{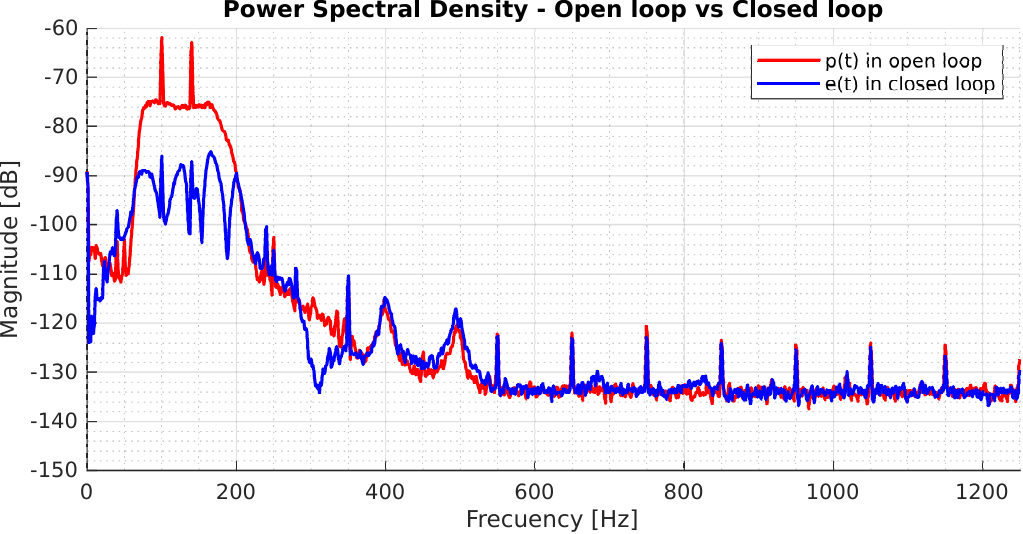} %
		\includegraphics[width=1\columnwidth]{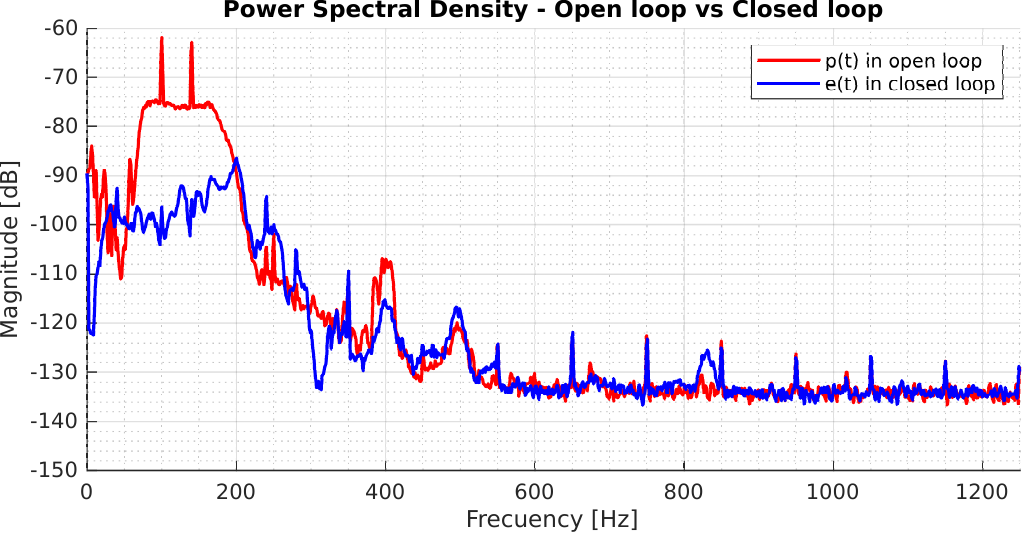} %
		
	\caption{Comparison of the residual noise PSD in open loop and in closed loop (t = 300 s) using the standard NLMS algorithm (top) and using the NLMS$+$ARIMA 2  algorithm (bottom). Disturbance: 70 - 170 Hz + 2 tonal disturbances (100 Hz, 140 Hz).}
	\label{fig:PSD_Comp}
\end{figure}

A signal consisting of a broad-band disturbance 70--170 Hz plus two tonal disturbances located at 100 Hz and 140 Hz was used as an unknown disturbance acting on the system. The steady state and transient attenuation\footnote{The attenuation is defined as the ratio between the variance of the residual noise in the absence of the control and the variance of the residual noise in the presence of the adaptive feedforward compensation. The variance is evaluated over an horizon of 3 seconds.} were evaluated for the various values of the parameters $c_1,~c_2$ and $d'_1$ given in Table~\ref{perf1}. The frequency characteristics of the various DAG given in Table~\ref{perf1} are shown in Fig. \ref{fig_bode_comp}. One observes that all are SPR (Phase lag between $0$ and $-90$ degrees with a $0$ dB average gain).
\begin{figure}[htb]
    \includegraphics[width=1\columnwidth]{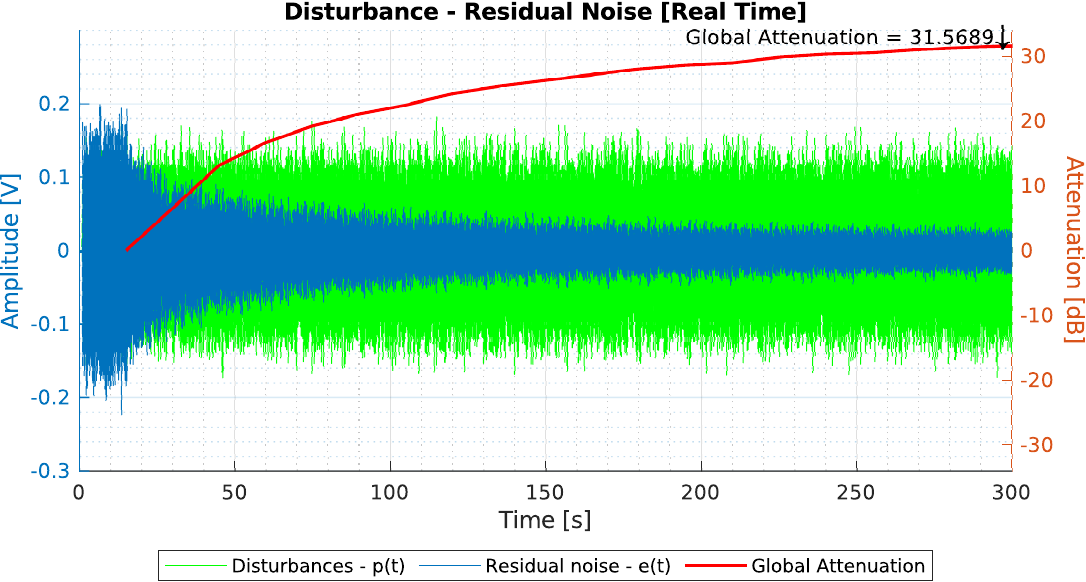}
		\includegraphics[width=1\columnwidth]{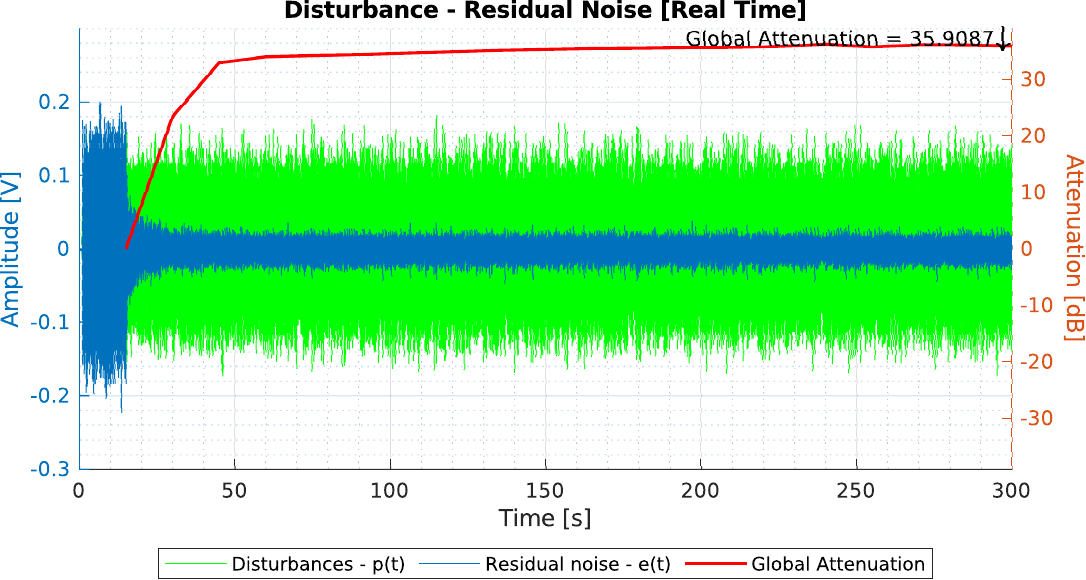}
	\caption{Time evolution of the residual noise and of the global attenuation using the standard NLMS adaptation algorithm (top) and using the NLMS + ARIMA2 algorithm (bottom), for $\mu=0.0003$.}
	\label{fig:resnoise}
\end{figure}
 The system was operated in open-loop during the first 15~s. Figure~\ref{fig:resnoise} shows the time response of the system as well as the evolution of the global attenuation when using the standard NLMS algorithm (top) and when the ARIMA 2 DAG is incorporated (bottom) with  $c_1=0.99,~c_2=0,~d'_1=0.9$ (last row of Table~\ref{perf1}). One observes a significant acceleration of the adaptation transient. 
  Figure~\ref{fig:PSD_Comp} shows for the same experiments a comparison of the PSD of the residual noise in the absence of the control and the PSD of the residual noise at 300 s under the control effect.  One can observes that the ARIMA 2 algorithm has improved the attenuation in the low frequency range with respect to the standard NLMS algorithm which is coherent with the frequency characteristic of the ARIMA 2 filter shown in Fig. \ref{fig_bode_comp}. 

\begin{figure}[htb]
    \includegraphics[width=1\columnwidth]{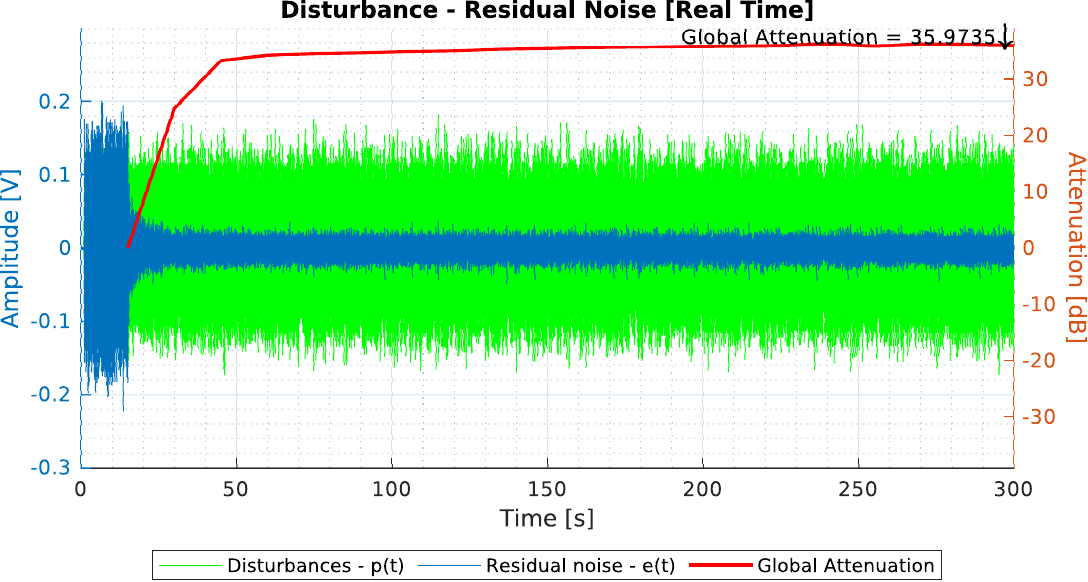}
		
	\caption{Time evolution of the residual noise and global attenutaion using the standard NLMS adaptation algorithm with $\mu=0.0075$.}
	\label{fig:HG_NLMS}
\end{figure}

  The acceleration obtained is equivalent to that obtained with an adaptation gain 25 times higher is used in the standard NLMS algorithms. This is illustrated in Figure~\ref{fig:HG_NLMS}. Similar results were obtained for the LMS and the PLMS algorithms.

\begin{figure}[!htb]
  \begin{center}
  \includegraphics[width=\columnwidth]{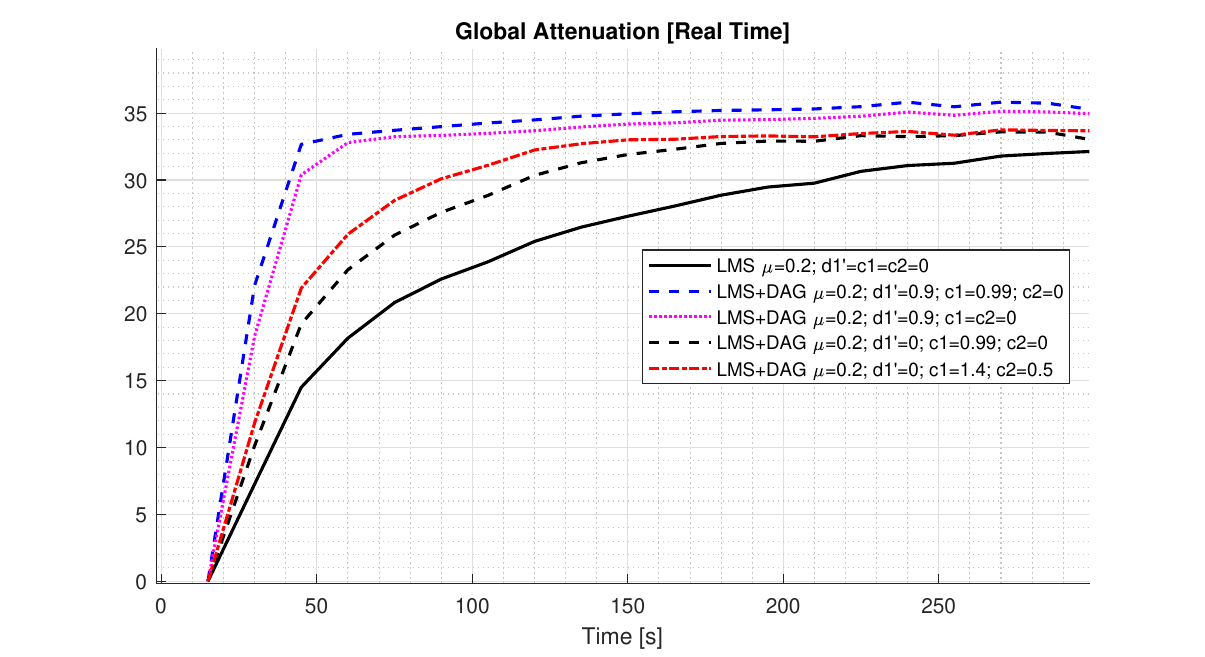}
  \caption{Time evolution of the global attenuation for the LMS algorithm with various DAGs, using $\mu=0.2$.}
  \label{fig_LMS_dag}
  \end{center}
\end{figure}

\begin{figure}[!htb]
  \begin{center}
  \includegraphics[width=\columnwidth]{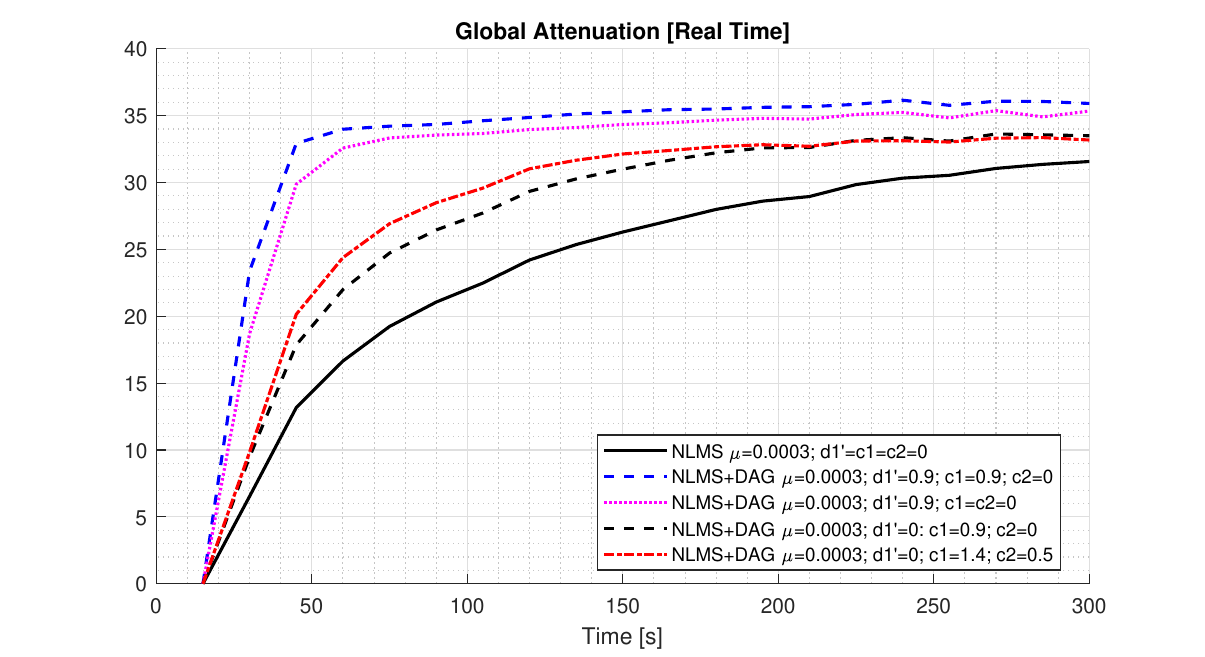}
  \caption{Time evolution of the global attenuation for the NLMS algorithm with various DAGs, using $\mu=0.0003$.}
  \label{fig_NLMS_dag}
  \end{center}
\end{figure}

\begin{figure}[!htb]
  \begin{center}
  \includegraphics[width=\columnwidth]{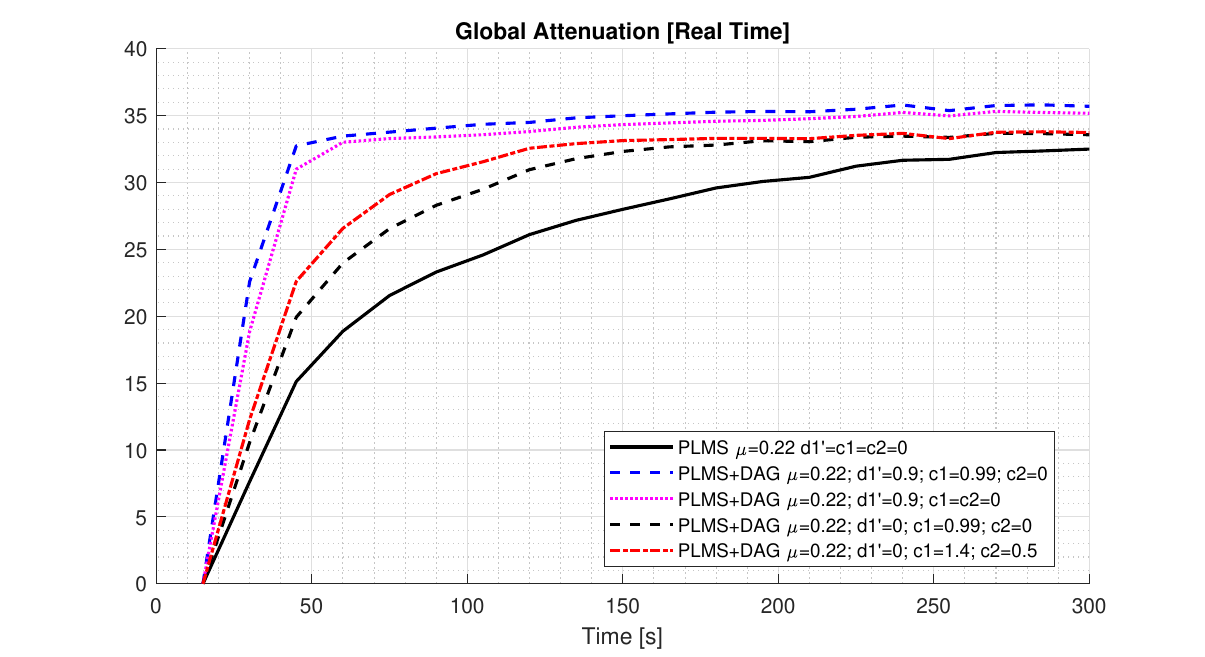}
  \caption{Time evolution of the global attenuation for the PLMS algorithm with various DAGs, using $\mu=0.22$.}
  \label{fig_PLMS_dag}
  \end{center}
\end{figure}

Figures \ref{fig_LMS_dag}, \ref{fig_NLMS_dag}, and \ref{fig_PLMS_dag} show the time evolution of the attenuation for the LMS, NLMS and PLMS algorithms when using the various DAG given in Table \ref{perf1}. As one can observe, the effect of the DAG is similar for the three algorithms. The DAG which has the highest steady state gain (ARIMA 2) gives the best result and the performance of the various DAGs are related to the value of their steady state gain (as predicted by Eq. \eqref{eq:rate}). 

Similar results were obtained for broad -- band disturbances with a wider spectrum. 
\section{Conclusion}
\label{Conclusions} 
The paper emphasizes the potential of the dynamic adaptation gain/step size (DAG) for improving the adaptation transients  of VS-LMS  adaptation algorithms.  The main point is that the DAG should be characterized by a \textit{strictly positive real} (SPR) transfer function
for a correct operation over the full frequency range from   $0$ to $0.5$ of the sampling frequency.
 Simulations results on  relevant adaptive filtering problems and experimental results on an adaptive active noise control system have illustrated the feasibility and the performance improvement achieved when a DAG is incorporated in VS-LMS adaptation algorithms.

\begin{appendices}
\section{Criterion for the selection of the coefficients of the ARIMA2 DAG}\label{Appen_SPR}
\begin{lem}\label{lemma_SPR}
The conditions assuring that $ H_{DAG}(z^{-1})=\frac{1+c_1z^{-1}+c_2z^{-2}}{1-d_1^{'}z^{-1}}$ is strictly positive real (SPR) are:
\begin{itemize}
\item for $c_2\leq 0$, $c_1$ must be such that
\[-1 -c_2 < c_1 < 1+c_2 \]

\item for  $c_2\geq 0$
\begin{itemize}
\item if the following condition is satisfied \[ 2(d_1^{'}-c_2) < \sqrt{2(c_2-c_2^2)(1-d_1^{'2})} <  2(d_1^{'}+c_2)  \]
  the maximum bound on $c_1$ is given by  \[c_1 < d_1^{'}-3d_1^{'}c_2+2\sqrt{2(c_2-c_2^2)(1-d_1^{'2})} \]
	otherwise the maximum bound on $c_1$ is given by \[ c_1 < 1+c_2 \]

\item   if the following condition is satisfied \[ 2(d_1^{'}-c_2) < -\sqrt{2(c_2-c_2^2)(1-d_1^{'2})} <  2(d_1^{'}+c_2)  \]
the minimum bound on $c_1$ is given by \[c_1 > d_1^{'}-3d_1^{'}c_2-\sqrt{2(c_2-c_2^2)(1-d_1^{'2})} \]
	otherwise the minimum bound on $c_1$ is given by \[ c_1 > -1-c_2 \]
\end{itemize}

\end{itemize}

\end{lem}
The proof of this result is given in \cite{LandauAuto23}.
\section{Proof of lemma \ref{lemstoch}} \label{proof_stoch1}
The verification of the Benveniste smoothness, boundness and mixing conditions
is similar to the verification done in \cite{Fan88} for the SHARF algorithm and it is omitted.
The focus is on the analysis of the Benveniste's ODE (ordinary differential equation) associated to the algorithm. The algorithm  of Eq. \eqref{eq:S2} can be rewritten as:
\begin{equation}
\hat{\mathbf{w}}(t+1)-  \hat{\mathbf{w}}(t)= \mu H_{DAG}(q^{-1})[\mathbf{r}(t,\hat{\mathbf{w}})e(t+1,\hat{\mathbf{w}})]
\label{aq:A1}
\end{equation}
i.e. one takes into account  that $\mathbf{r}(t)$ may depend upon $\hat{\mathbf{w}}$ and that clearly $e(t+1)$ depends upon $\hat{\mathbf{w}}$.
The associated ODE equation is:
\begin{equation}
\frac{d\hat{\mathbf{w}}}{d\tau}=\mu \mathbf{E}\lbrace H_{DAG}[\mathbf{r}(t,\hat{\mathbf{w})}e(t+1,\hat{\mathbf{w}})]\rbrace
\label{aq:A2}
\end{equation}
But from Eqs.  \eqref{eq:stoch 1} and \eqref{eq_aposteriori} one gets:
\begin{equation}
e(t+1,\hat{\mathbf{w}})=\mathbf{r}^T(\mathbf{w} -\hat{\mathbf{w}}) +\mathbf{n}(t+1)
\label{aq:A3}
\end{equation}
which leads to:
\begin{align}
\frac{d\hat{\mathbf{w}}}{dt}&=-\mu \mathbf{E}\lbrace H_{DAG}(q^{-1})[\mathbf{r}(n,\hat{\mathbf{w}}) \mathbf{r}(n,\hat{\mathbf{w}})^T]\rbrace(\hat{\mathbf{w}}-\mathbf{w})\nonumber \\
&+ \mu \mathbf{E}\lbrace H_{DAG}(q^{-1})[\mathbf{r}(n,\hat{\mathbf{w}}) \mathbf{n}(t+1)]\rbrace
\label{eq:A4}
\end{align}
Since $\mathbf{n}(t+1)$ is a white noise or is a stochastic process independent of $\mathbf{r}(t,\hat{\mathbf{w}})$, the second term of the right hand side of Eq. \eqref{eq:A4} will be null. Therefore, the properties of the algorithm for large $\tau$ and small $\mu$ will be determined by the stability of:
\begin{equation}
\frac{d\hat{\mathbf{w}}}{d\tau}=-\mu G_{\mathbf{w}}(\hat{\mathbf{w}}-\mathbf{w})
\label{eq:A5}
\end{equation} 
where:
\begin{align}
 G_{\mathbf{w}}&=\mathbf{E}\lbrace H_{DAG}(q^{-1})[\mathbf{r}(n,\hat{\mathbf{w}}), \mathbf{r}(n,\hat{\mathbf{w}})^T]\rbrace\\
&=H_{DAG}(q^{-1})\mathbf {E}_{\mathbf{r}}
 \label{eq:Gw}
\end{align}
where $\mathbf {E}_{\mathbf{r}}= \mathbf{E}[\mathbf{r}(n,\hat{\mathbf{w}}), \mathbf{r}(n,\hat{\mathbf{w}})^T]$.
 The asymptotic stability of the above equation can be assessed using a Lyapunov function candidate: 
\begin{equation}
V(\hat{\mathbf{w}})=(\hat{\mathbf{w}}-\mathbf{w})^T\mu^{-1}(\hat{\mathbf{w}}-\mathbf{w})
\label{eq:V}
\end{equation}
The derivative of the Lyapunov function candidate along the trajectories of the system has the form:
\begin{equation}
\dot{V}(\hat{\mathbf{w}})=-(\hat{\mathbf{w}}-\mathbf{w})^T [G_{\mathbf{w}}+G_{\mathbf{w}}^T](\hat{\mathbf{w}}-\mathbf{w})
\label{eq:Vdot}
\end{equation}
It remains to show that $G_{\mathbf{w}}+G_{\mathbf{w}}^T$ is positive definite. To do this, we will use the integration on the unit circle. One obtains:
\begin{align}
G_{\mathbf{w}}+ G_{\mathbf{w}}=& \frac{1}{2\pi}\int_{-\pi}^{\pi}[H_{DAG}(e^{j\omega})\mathbf{\Psi}(\omega)\\
&+\mathbf{\Psi}(\omega)^TH_{DAG}(e^{-j\omega})]d\omega\\
=& \frac{1}{2\pi}\int_{-\pi}^{\pi}2 Re[H^{ii}_{DAG}(e^{\omega})]\mathbf{\mathbf{\Psi}}(\omega)d\omega 
\end{align}
taking into account that $\mathbf{\Psi}(\omega)$, the power spectral density of $\mathbf{r}(t,\mathbf{w}$) is symmetric. Since $H^{ii}_{DAG}$ is SPR and $\mathbf{\Psi}(\omega)$ is positive definite (as a consequence of condition of Eq. \eqref{eq:S4}), one concludes that $G_{\mathbf{w}}+G_{\mathbf{w}}^T$ is positive definite.

\end{appendices}

\bibliographystyle{IEEEtranS}        %
\bibliography{Bibliography}      %

\end{document}